\documentclass[10pt,a4paper,leqno]{amsart}
\usepackage[english]{babel}
\usepackage[T1]{fontenc}
\usepackage{amssymb}
\usepackage{amsfonts}
\usepackage{amsmath}
\usepackage[T1]{fontenc}
\usepackage{eucal}

\begin{document}

%
%
\null{}

\vspace{5ex}

\title{\Large Subharmonic functions, generalizations and separately  subharmonic functions}
\author{Juhani Riihentaus}
\date{}
\maketitle
\vspace{-0.35in}
\begin{center}
{Department of Physics and Mathematics, University of Joensuu\\
P.O. Box 111, FI-80101 Joensuu, Finland \\
juhani.riihentaus@joensuu.fi}
\end{center}
\vspace*{4ex}

\noindent{\emph{Abstract:}} First, we give the definition for quasi-nearly subharmonic functions, now for general, not necessarily nonnegative functions, 
unlike previously. We point out that our function class incudes, among others, 
quasisubharmonic functions, nearly subharmonic functions (in a slightly generalized sense) and almost subharmonic functions. We also give some 
basic properties of quasi-nearly subharmonic functions. Second, after recalling some  of the existing subharmonicity  results of separately 
subharmonic functions,  we  give the corresponding counterparts for  separately quasi-nearly subharmonic functions, thus improving  previous results of 
ours, of Lelong, of Avanissian and of Arsove. Third,  we give two  results concerning the subharmonicity of a function  subharmonic 
with respect to the first variable and harmonic with respect to the second variable. The first  result improves a result of Arsove,  
concerning the case when the function has, in addition, locally a negative integrable minorant. The second result improves a result of Ko\l odziej and 
Thorbi\"ornson concerning the subharmonicity of a function 
subharmonic and ${\mathcal{C}}^2$ in the first variable and harmonic in the second.

\vspace{0.5ex}


\vspace{4ex}

\noindent{{\textbf{1.\,\, Introduction}}}

\vspace{3ex}

\noindent {\textbf{1.1. Separately subharmonic functions.}} Solving a long standing problem, Wiegerinck [Wi88, Theorem, p.~770], see also [WZ91, Theorem~1, p.~246], 
   showed that a separately subharmonic function need not
be subharmonic. On the other hand, Armitage and Gardiner [AG93, Theorem~1, p.~256] showed that a separately subharmonic function $u$ in a domain $\Omega $
in ${\mathbb{R}}^{m+n}$, $m\geq n\geq 2$, is subharmonic provided 
$\phi (\log^+u^+)$ is locally integrable in $\Omega$, where $\phi : [0,+\infty )\rightarrow [0,+\infty )$ is an increasing function such that 
\begin{equation}\int\limits_1^{+\infty }s^{(n-1)/(m-1)}(\phi (s))^{-1/(m-1)}\,ds<+\infty .\end{equation}
Armitage's and Gardiner's result includes the previous  results of Lelong [Le45,  Théorème~1~bis, p.~315], of Avanissian 
[Av61, Théorème~9, p.~140], see also [Her71, Theorem, p.~31], of  Arsove [Ar66, Theorem~1, p.~622] and ours:

\vspace{1ex}

\noindent{\textbf{Theorem A.}} ([Ri89, Theorem~1, p. 69]) \emph{ Let $\Omega $ be a domain in ${\mathbb{R}}^{m+n}$, \mbox{$m,n\geq 2$.} 
Let $u:\, \Omega \rightarrow 
[-\infty ,+\infty )$ be such that} 
 \begin{itemize}
\item[(a)] \emph{for each $y\in {\mathbb{R}}^n$ the function} 
\[\Omega (y)\ni x\mapsto u(x,y)\in [-\infty ,+\infty )\]
\emph{is subharmonic,}
 \item[(b)] \emph{for each $x\in {\mathbb{R}}^m$ the function} 
\[\Omega (x)\ni y\mapsto u(x,y)\in [-\infty ,+\infty )\]
\emph{is subharmonic,}
\item[(c)] \emph{for some $p>0$ there is a function}   $v\in {\mathcal{L}}_{\textrm{loc}}^p(\Omega )$ \emph{such that $u\leq v$}.
\end{itemize} 
 \emph{Then} $u$ 
\emph{is subharmonic in $\Omega $}.

\vspace{1ex}

Instead of (c), Lelong and Avanissian used the stronger assumption that the function is locally bounded above, and Arsove the assumption that  the 
function has a locally  integrable majorant. Though the result of Armitage and Gardiner is even \,``almost''  sharp, it is, nevertheless,  relevant to try to 
find both new proofs and generalizations   also to the previous results. This is justified because of two reasons. First, 
our ${\mathcal{L}}^p_{\textrm{loc}}$-integrability condition, \mbox{$p>0$,}  is, unlike the condition of 
Armitage and Gardiner  (1), very simple, and second, more importantly, Armitage's and Gardiner's proof is \emph{based}  on the results of Lelong and Avanissian,
or, alternatively, on the  more general result of 
Arsove or on  our Theorem~A above, say.

\vspace{2ex}

\noindent\textbf{1.2. Functions subharmonic in one variable and harmonic in the other.} An  open problem is, whether a function which is subharmonic in
 one variable and harmonic in the other, is subharmonic. For results on this area, 
see e.g. [Ar66], [Im90], [WZ91], [CS93] and [KT96] and the references therein. We consider here a result of Arsove, Theorem~B below, 
and a  result of Ko\l odziej and 
Thorbi\"ornson, Theorem~C below. Observe that the situation here is similar with that pointed out above in  {\textbf{1.1}}:  The proofs of these theorems (and also the 
proofs of our improvements, see Theorems~4 and 5 below),   are  again based on the cited result of Lelong and Avanissian, or, alternatively, on the results 
of Arsove and ours, say.  

\vspace{1ex}

\noindent{\textbf{Theorem B.}} ([Ar66, Theorem~2, p.~622]) \emph{ Let $\Omega $ be a domain in ${\mathbb{R}}^{m+n}$, \mbox{$m,n\geq 2$.} 
Let $u:\, \Omega \rightarrow 
{\mathbb{R}}$ be such that} 
 \begin{itemize}
\item[(a)] \emph{for each $y\in {\mathbb{R}}^n$ the function} 
\[\Omega (y)\ni x\mapsto u(x,y)\in {\mathbb{R}}\]
\emph{is subharmonic,}
 \item[(b)] \emph{for each $x\in {\mathbb{R}}^m$ the function} 
\[\Omega (x)\ni y\mapsto u(x,y)\in {\mathbb{R}}\]
\emph{is harmonic,}
\item[(c)] \emph{there is a nonnegative function } $\varphi \in {\mathcal{L}}^1_{\textrm{loc}}(\Omega )$ \emph{such that} $-\varphi \leq u$.
\end{itemize}
\emph{Then} $u$ 
\emph{is subharmonic in $\Omega $}.
\vspace{1ex}

Arsove's proof was based on mean value operators. Much later  Cegrell and Sadullaev 
[CS93, Theorem~3.1, p. 82] gave a new proof using Poisson modification.  In Theorem~2, in its Corollary  and in Theorem~3 below we give  
 generalizations to this result.

\vspace{1ex}

Ko\l odziej and Thorbi$\ddot {\textrm{o}}$rnson gave the following result. 

\vspace{1ex}

\noindent{\textbf{Theorem~C.}} ([KT96, Theorem~1, p.~463]) \emph{ Let $\Omega $ be a domain in ${\mathbb{R}}^{m+n}$, \mbox{$m,n\geq 2$.} 
Let $u:\, \Omega \rightarrow 
{\mathbb{R}}$ be such that} 
 \begin{itemize}
\item[(a)] \emph{for each $y\in {\mathbb{R}}^n$ the function} 
\[\Omega (y)\ni x\mapsto u(x,y)\in {\mathbb{R}}\]
\emph{is subharmonic and ${\mathcal{C}}^2$,}
 \item[(b)] \emph{for each $x\in {\mathbb{R}}^m$ the function} 
\[\Omega (x)\ni y\mapsto u(x,y)\in {\mathbb{R}}\]
\emph{is harmonic.}
\end{itemize}
\emph{Then} $u$ 
\emph{is subharmonic and continuous in $\Omega $}.

\vspace{1ex}

In [Ri07$_1$, Theorem~6, p.~233] and [Ri07$_2$, Theorem~1, p.~438, and Theorem~2, pp.~442-443] we improved the above result of Ko\l odziej and 
Thorbi\"ornson.
Now  in Theorem~4 and  Theorem~5 below we improve these results still further. Instead of the standard Laplacians of 
$\mathcal{C}^2$ functions we will use  generalized Laplacians, that is, the Blaschke-Privalov operators.

\vspace{4ex}

\noindent{\textbf{2. \,\, Notation, subharmonic functions, generalizations,  and properties}} 

\vspace{3ex}

\noindent {\textbf{2.1. Notation.}} Our notation is rather standard, see e.g. [Ri06$_1$] and [Her71]. $m_N$ is the Lebesgue measure 
in the Euclidean space ${\mathbb{R}}^N$, $N\geq 2$. We write $\nu _N$ for the Lebesgue measure of the unit ball $B^N(0,1)$ 
in ${\mathbb{R}}^N$, thus $\nu _N=m_N(B^N(0,1))$. $D$ is a domain of ${\mathbb{R}}^N$.
The complex space ${\mathbb{C}}^n$ is identified with the real space  ${\mathbb{R}}^{2n}$, $n\geq 1$.  
Constants will be denoted by $C$ and $K$. They will be nonnegative and may vary from line to line.

\vspace{1ex}

\noindent {\textbf{2.2. Subharmonic functions and generalizations.}} We recall that an upper semicontinuous function $u:\, D\rightarrow [-\infty ,+\infty )$ is \emph{subharmonic} if 
for all $\overline{B^N(x,r)}\subset D$,
\[u(x)\leq \frac{1}{\nu _N\, r^N}\int\limits_{B^N(x,r)}u(y)\, dm_N(y).\]
The function $u\equiv -\infty $  is considered  subharmonic. 

We say that a function 
$u:\, D\rightarrow [-\infty ,+\infty )$ is \emph{nearly subharmonic}, if $u$ is Lebesgue measurable, $u^+\in {\mathcal{L}}^1_{\textrm{loc}}(D)$, 
and for all $\overline{B^N(x,r)}\subset D$,  
\begin{equation}u(x)\leq \frac{1}{\nu _N\, r^N}\int\limits_{B^N(x,r)}u(y)\, dm_N(y).\end{equation}
Observe that in the standard definition of nearly subharmonic functions one uses the slightly stronger  assumption that 
$u\in {\mathcal{L}}^1_{\textrm{loc}}(D)$, see e.g. [Her71, p.~14]. However, our above, slightly 
more general definition seems to be  more practical, see  below  Proposition~1~(iii) and Proposition~2~(vi) and (vii).
Proceeding as in [Her71, proof of Theorem~1, pp.~14-15] (and referring also to Proposition~1~(iii) and Proposition~2~(vii), see below) one obtains the following result:
 
\vspace{1ex}

\noindent{\textbf{Lemma.}} \emph{Let $D$ be a domain in ${\mathbb{R}}^N$, $N\geq 2$. Let $u:\,D\rightarrow [-\infty ,+\infty )$ be Lebesgue measurable. 
Then  $u$ is nearly subharmonic in $D$ if and only if there exists a function $u^*$, subharmonic in $D$  such that $u^*\geq u$ and $u^*=u$ almost 
everywhere in $D$. Here $u^*$ is the lowest  upper semicontinuous majorant of $u$:
\begin{displaymath}u^*(x)=\limsup_{x'\rightarrow x}u(x').\end{displaymath}
$u^*$ is called the regularized subharmonic function to $u$.}

\vspace{1ex} 

Observe also that \emph{almost subharmonic functions,}
in the sense of Szpilrajn [Sz33] (see also [Ra37, 3.30, p.~20] and  [LL01, p.~238]; Lieb and Loss  even  call this class briefly subharmonic functions!),
are, more or less,  included in our definition of nearly subharmonic functions, in the following sense. Let  $u: \, D\rightarrow [-\infty ,+\infty )$ be 
almost subharmonic, that is $u\in {\mathcal{L}}^1_{\textrm{loc}}(D)$ and for almost every $x\in D$ and for every 
$r>0$ such that $\overline{B(x,r)}\subset D$ the  mean value inequality (2) is satisfied. Let
\[D_1:=\{x\in D\,:\, u(x)\leq \frac{1}{\nu _N\, r^N}\int\limits_{B^N(x,r)}u(x')\, dm_N(x')\, {\textrm{ for all }}\overline{B^N(x,r)}\subset D\,   \}\]
and let $A:=D\setminus D_1$. Define $\tilde {u}:\, D\rightarrow [-\infty ,+\infty )$,
\[\tilde {u}(x):=\begin{cases}u(x) {\textrm{ when}}&x\in D_1,\\
-\infty {\textrm{ when}}&x\in A.\end{cases}\]
Since by assumption $m_N(A)=0$, it is easy to see that $\tilde {u}$ is nearly subharmonic in $D$.

\vspace{1ex}

The previous definition for quasi-nearly subharmonic functions was restricted to nonnegative functions, see [Pa94], [Mi96], [Ri00], [PR08], [Ri05], [Ri06$_1$].
Now we say that a Lebesgue measurable function $u:\,D \rightarrow 
[-\infty ,+\infty )$ is \emph{$K$-quasi-nearly subharmonic}, if  $u^+\in{\mathcal{L}}^{1}_{\textrm{loc}}(D)$ and if there is a 
constant $K=K(N,u,D)\geq 1$
such that for all $x\in D$ and $r>0$ such that  $\overline{B^N(x,r)}\subset D$, one has   
\begin{equation*} u_M(x)\leq \frac{K}{\nu _N\,r^N}\int\limits_{B^N(x,r)}u_M(y)\, dm_N(y)\end{equation*}
for all $M\geq 0$, where $u_M:=\max\{u,-M\}+M$. A function $u:\, D\rightarrow [-\infty ,+\infty )$ is \emph{quasi-nearly subharmonic}, if $u$ is 
$K$-quasi-nearly subharmonic in $D$ for some $K\geq 1$.

\vspace{1ex}

In addition to the above defined class of quasi-nearly subharmonic functions, we will  consider also a proper subclass.  A Lebesgue  measurable function 
$u:\,D \rightarrow [-\infty ,+\infty )$ is \emph{$K$-quasi-nearly subharmonic n.s. (in the narrow sense)}, if $u^+\in{\mathcal{L}}^{1}_{\textrm{loc}}(D)$ and if there is a 
constant $K=K(N,u,D)\geq 1$
such that for all $\overline{B^N(x,r)}\subset D$, one has   
\begin{equation*} u(x)\leq \frac{K}{\nu _N\,r^N}\int\limits_{B^N(x,r)}u(y)\, dm_N(y).\end{equation*}
A function $u:\, D\rightarrow [-\infty ,+\infty )$ is \emph{quasi-nearly subharmonic n.s.}, if $u$ is 
$K$-quasi-nearly subharmonic n.s. in $D$ for some $K\geq 1$. 

Observe that already Domar in [Do57, p.~430] has pointed out the relevance of the class of (nonnegative) quasi-nearly subharmonic functions. 
For, at least partly,  an even more general function class, see [Do88].

\vspace{1ex}

As an example of a subclass of quasi-nearly subharmonic functions,  and as a counterpart to nonnegative harmonic functions we recall the definition of Harnack functions, see [Vu82, p.~259]. 
Suppose for a while that $D\ne {\mathbb{R}}^N$. A continuous function 
$u:\, D\rightarrow [0,+\infty )$ is a $\lambda $-\emph{Harnack function}, if there are constants $\lambda \in (0,1)$ and $C_\lambda =C(\lambda )\geq 1$ such that 
\begin{equation*}\max_{z\in \overline{B^N(x,\lambda r)}}u(z)\leq C_\lambda \, \min_{z\in \overline{B^N(x,\lambda r)}}u(z)\end{equation*}
whenever $\overline{B^N(x,r)}\subset D$. It is well-known that for each compact set $F\ne \emptyset$ in $D$ there exists a smallest constant 
$C(F)\geq C_\lambda $ depending only on $N$, 
$\lambda $, $C_\lambda $ and $F$ such that for all $u$ satisfying the above condition,
\begin{equation*}\max_{z\in F}u(z)\leq C(F)\, \min_{z\in F}u(z).\end{equation*}
A continuous function $u:\, D\rightarrow [0,+\infty )$ is a \emph{Harnack function}, if it is a $\lambda $-Harnack function for some $\lambda \in (0,1)$.

\vspace{2ex}

\noindent{\textbf{2.3.}} The defined function  classes,  quasi-nearly subharmonic functions,  quasi-nearly subharmonic functions n.s. and 
Harnack functions, are all natural,  they have important and interesting properties and, at the same time, they are large,      
See e.g. [Pa94], [Mi96], [Ri00], [PR08], [Ri05], [Ri06$_1$], [Vu82], and Propositions~1 and~2 below. We recall here only that these function classes 
include, among others,  subharmonic functions, and, more generally,  quasisubharmonic (see e.g. [Br38], [Le45, p.~309], [Av61, p.~136], [Her71, p.~26]) and also 
nearly subharmonic functions (see e.g. [Her71, p.~14]),    also functions satisfying certain natural  growth conditions, especially  
certain eigenfunctions, and  polyharmonic functions. Especially, the class of Harnack functions includes, among others, nonnegative harmonic functions 
as well as nonnegative solutions of some elliptic equations. In particular, the partial differential equations associated with quasiregular mappings 
belong to this family of elliptic equations, see [Vu82].

To motivate the above defined function classes, the class  of quasi-nearly subharmonic functions and the class of quasi-nearly subharmonic n.s. functions, 
even more, we give below four simple examples.

\vspace{1ex}

\noindent{\textbf{Example~1.}} Let $D$ be a domain in ${\mathbb{R}}^N$, $N\geq 2$. Any Lebesgue measurable function   $u:\,D\rightarrow [m,M]$, 
where $0<m\leq M<+\infty $, is quasi-nearly subharmonic, and, because of Proposition~1~(i), also quasi-nearly subharmonic n.s. If $u$ is moreover 
continuous, it is even a Harnack function, see [Vu82, pp.~259, 263].

\vspace{1ex}

\noindent{\textbf{Example~2.}} The function  $u:\,{\mathbb{R}}^2\rightarrow {\mathbb{R}}$,
\begin{displaymath} u(x,y):=\begin{cases}-1, & {\textrm{when }}\, y<0,\\
1, & {\textrm{when }}\, y\geq 0,\end{cases}\end{displaymath}
is $2$-quasi-nearly subharmonic, but not quasi-nearly subharmonic n.s. 

\vspace{1ex}

\noindent{\textbf{Example~3.}} Let $D=(0,2)\times (0,1)$, let $c<0$ be  arbitrary. Let $E\subset D$ be a Borel set of zero Lebesgue measure. Let 
$u:\,D\rightarrow [-\infty ,+\infty )$,
\begin{displaymath} u(x,y):=\begin{cases}c , & {\textrm{ when }}\, (x,y)\in E,\\
1, &{\textrm{ when }}\, (x,y)\in D\setminus E  {\textrm{ and }}\, 0<x<1,\\
2, & {\textrm{ when }}\, (x,y)\in D\setminus E {\textrm{ and }}\, 1\leq x<2.\end{cases}\end{displaymath}
The function $u$ attains both negative and positive values, it is $2$-quasi-nearly subharmonic n.s, but not nearly subharmonic. 
Recall  that  previously we have considered only nonnegative quasi-nearly subharmonic functions. 

\vspace{1ex}

\noindent{\textbf{Example~4.}} Let $D$ be a domain in ${\mathbb{R}}^N$, $N\geq 2$, and let $u:\,D\rightarrow [-\infty ,+\infty )$ be any quasi-nearly 
subharmonic function n.s. Let $E\subset D$ be a Borel set of zero Lebesgue measure. Let 
$v:\,D\rightarrow [-\infty ,+\infty )$,
\begin{displaymath} v(x,y):=\begin{cases}-\infty , & {\textrm{ when }}\, (x,y)\in E,\\
u(x,y), &{\textrm{ when }}\, (x,y)\in D\setminus E.  \end{cases}\end{displaymath}
The function $v$ is clearly quasi-nearly subharmonic n.s. 

\vspace{1ex}

For the convenience of the reader we recall the following definition, see also [PR08, Lemma~1 and Remark~1, p.~93].

A function $\psi :\, [0,+\infty )\rightarrow [0,+\infty )$ is {\textit {permissible}}, if  
there  exists an increasing (strictly or not), convex function
$\psi_1 :\,[0,+\infty )\rightarrow [0,+\infty )$
and a strictly  increasing surjection
\mbox{$\psi_2 :\, [0,+\infty )\rightarrow [0,+\infty )$} such that $\psi =\psi_2\circ \psi_1$ and such that the following
conditions are satisfied:
\begin{itemize}
\item[(a)] $\psi_1$  satisfies the
$\varDelta_2$-condition,  
\item[(b)] $\psi_2^{-1}$ satisfies the
$\varDelta_2$-condition,  
\item[(c)] the function $t\mapsto \frac{\psi_2(t)}{t}$ is
{\textit {quasi-decreasing}}, i.e. there is a constant $C=C(\psi_2)>0$ such
that
\begin{displaymath}\frac{\psi_2(s)}{s}\geq  C\,    \frac{\psi_2(t)}{t}\end{displaymath}
for all $0\leq s\leq t$.
\end{itemize}
\noindent Recall that a function $\varphi :\,[0,+\infty )\rightarrow [0,+\infty )$ satisfies the $\varDelta_2$-\emph{condition}, if there is a constant $C=C(\varphi )\geq 1$
such that $\varphi (2t)\leq C\,\varphi (t)$ for all $t\in [0,+\infty )$. 

\vspace{1ex}

Examples of permissible functions are:  $\psi_1(t)=t^p$, $p>0$, and   $\psi_2(t)=c\, t^{p\alpha }[\log (\delta +t^{p\gamma })]^\beta $, $c>0$, 
$0<\alpha <1$, $\delta \geq 1$, $\beta ,\gamma \in {\mathbb{R}}$ such that $0<\alpha +\beta \,\gamma <1$, and $p\geq 1$. And also functions of the form 
$\psi _3=\phi \circ \varphi$, where $\phi :[0,+\infty )\rightarrow [0,+\infty )$ is a concave surjection whose inverse $\phi ^{-1}$ satisfies 
the $\Delta_2$-condition and $\varphi :[0,+\infty )\rightarrow [0,+\infty )$ is an increasing, convex function satisfying the $\Delta_2$-condition.

\noindent{\textbf{2.4. Properties.}} In Proposition~1 and Proposition~2 below we collect some of the basic properties of the above defined function classes. 
For the properties (iv) of Proposition~1, and (iii), (v) and (vi) of Proposition~2, see already  [PR08, Proposition~1, Theorem~A, Corollary~1 and Theorem~B, p.~91],
[Ri89, Lemma, p.~69] and the references therein.   

\vspace{1ex}

\noindent{\textbf{Proposition~1.}} \emph{Let $D$ be a domain in ${\mathbb{R}}^N$, $N\geq 2$.}
\begin{itemize}
\item[(i)] \emph{If  $u:\, D\rightarrow [0 ,+\infty )$   is Lebesgue measurable and} $u\in {\mathcal{L}}^1_{\textrm{loc}}(D)$\emph{, then $u$ is
 $K$-quasi-nearly subharmonic if and only if $u$ is $K$-quasi-nearly subharmonic n.s., that is, if 
for all $\overline{B^N(x,r)}\subset D$,}
\[u(x)\leq \frac{K}{\nu _N\, r^N}\int\limits_{B^N(x,r)}u(y)\, dm_N(y).\]
\item[(ii)] \emph{If  $u:\, D\rightarrow [-\infty  ,+\infty )$   is $K$-quasi-nearly subharmonic n.s., then $u$ is $K$-quasi-nearly subharmonic in $D$, 
but not necessarily conversely.}
\item[(iii)] \emph{A function $u:\, D\rightarrow [-\infty ,+\infty )$ is $1$-quasi-nearly subharmonic if and 
only if it is nearly subharmonic, that is, it is  $1$-quasi-nearly subharmonic n.s.}
\item[(iv)] \emph{If $u:\, D\rightarrow [0,+\infty )$ is quasi-nearly subharmonic and $\psi :\, [0,+\infty )\rightarrow [0,+\infty )$ is permissible, 
then $\psi \circ u$ is quasi-nearly subharmonic in $D$.} 
\item[(v)] \emph{Harnack functions are quasi-nearly subharmonic.}
\end{itemize} 

\vspace{1ex}

\noindent\emph{Proof.} For the proof of (iv), see [Ri$06_1$, Lemma~2.1, p.~32]. We leave the cases (i), (ii), (v) to the reader. To prove the case (iii) 
suppose that  $u$ is nearly subharmonic in $D$. Then clearly $u_M$ is nearly subharmonic for all $M\geq 0$, and 
thus for all $x\in D$ and $r>0$ such that  $\overline{B^N(x,r)}\subset D$, one has   
\begin{equation*} u_M(x)\leq \frac{1}{\nu _N\,r^N}\int\limits_{B(x,r)}u_M(y)\, dm_N(y).\end{equation*}
Hence $u$ is $1$-quasi-nearly subharmonic.

On the other hand, if $u$ is $1$-quasi-nearly subharmonic in $D$, then one sees at once, with the aid of Lebesgue Monotone Convergence Theorem, that
$u$ is nearly subharmonic in $D$. \hfill \qed

\vspace{1ex}

\noindent{\textbf{Proposition~2.}} \emph{Let $D$ be a domain in ${\mathbb{R}}^N$, $N\geq 2$.}
\begin{itemize}
 \item[(i)] \emph{If  $u:\, D\rightarrow [-\infty ,+\infty )$   is $K_1$-quasi-nearly subharmonic and $K_2\geq K_1$, then 
$u$ is $K_2$-quasi-nearly subharmonic in $D$.}
\item[(ii)] \emph{If $u_1:\, D\rightarrow [-\infty ,+\infty )$    and $u_2:\, D\rightarrow [-\infty ,+\infty )$   are $K$-quasi-nearly subharmonic n.s.,
  then $\lambda_1 u_1+\lambda _2u_2$ is $K$-quasi-nearly subharmonic n.s. in $D$ for all  $\lambda_1,\, \lambda_2   \geq 0$.} 
\item[(iii)] \emph{If  $u:\, D\rightarrow [-\infty ,+\infty )$ is quasi-nearly subharmonic, then $u$ is locally bounded above in $D$.}
\item[(iv)] \emph{If $u_j:\, D\rightarrow [-\infty ,+\infty )$, $j=1,2,\dots$, are $K$-quasi-nearly subharmonic (resp. $K$-quasi-nearly subharmonic n.s.), 
and $u_j\searrow u$ as $j\rightarrow +\infty $,
then $u$ is $K$-quasi-nearly subharmonic (resp. $K$-quasi-nearly subharmonic n.s.) in $D$.}
\item[(v)] \emph{If $u:\, D\rightarrow [-\infty ,+\infty )$   is $K_1$-quasi-nearly subharmonic and 
$v:\, D\rightarrow [-\infty ,+\infty )$ is $K_2$-quasi-nearly subharmonic, then $\max\{u,v\}$ is $\max\{K_1,K_2\}$-quasi-nearly 
subharmonic in $D$. Especially,  
 $u^+:=\max\{u,0\}$ is $K_1$-quasi-nearly subharmonic in $D$.}
\item[(vi)] \emph{Let ${\mathcal{F}}$ be a family of  $K$-quasi-nearly subharmonic (resp. $K$-quasi-nearly subharmonic n.s.) functions in $D$ and let 
$w:=\sup_{u\in {\mathcal{F}}}u$. If $w$ is Lebesgue measurable and} $w^+\in {\mathcal{L}}_{{\textrm{loc}}}^1(D)$, \emph{then $w$ is 
$K$-quasi-nearly subharmonic (resp. $K$-quasi-nearly subharmonic n.s.) in $D$.}
\item[(vii)] \emph{If  $u:\, D\rightarrow [-\infty ,+\infty )$   is quasi-nearly subharmonic n.s., then either $u\equiv -\infty $ or $u$ is finite almost 
everywhere in $D$, and} 
$u\in {\mathcal{L}}^1_{\textrm{loc}}(D)$.
\end{itemize}

\vspace{1ex}

\noindent{\textbf{Remark.}} Related to (ii) above, it is easy to see that, if $u:\, D\rightarrow [-\infty ,+\infty )$  is $K$-quasi-nearly subharmonic,
  then $\lambda u+C$ is $K$-quasi-nearly subharmonic in $D$ for all  $\lambda  \geq 0$ and $C\geq 0$. However, at present we do not know 
whether $\lambda_1 u_1+\lambda _2u_2$ is $K$-quasi-nearly subharmonic provided 
$u_1:\, D\rightarrow [-\infty ,+\infty )$    and $u_2:\, D\rightarrow [-\infty ,+\infty )$   are $K$-quasi-nearly subharmonic 
 and  $\lambda_1,\, \lambda_2 \geq 0$. 

\vspace{1ex}

\noindent\emph{Proof of Proposition~2.} We leave the rather easy cases (i)--(vi) to the reader and prove only (vii). Our proof is nearly verbatim the same 
as [Hel69, proof of Theorem~4.10, p. 66], except perhaps in a couple of the last lines of the proof.

By (vi) $u^+$ is quasi-nearly subharmonic n.s., and by (iii) locally bounded above in $D$. Therefore, for each  $\overline{B(x,r)}\subset D$, the integral
\[\int\limits_{B(x,r)}u(y)dm_N(y)\]
is defined either as $-\infty$ or as a finite real number. Suppose that $-\infty <u(x_0)$ for some $x_0\in D$. Then
\begin{displaymath}-\infty <u(x_0)\leq \frac{K}{\nu _Nr^N}\int\limits_{B^N(x_0,r)}u(y)\, dm_N(y)<+\infty .\end{displaymath}
Therefore there exists $r_0>0$ such that $u(x)\in {\mathbb{R}}$ for almost every $x\in B^N(x_0,r_0)$.  

Write 
\[H:=\{\, x\in D\,:\, u(y) {\textrm{ finite for almost every }}y\in B(x,r)\,{\textrm{ for some }}r>0\,\textrm{ s.t. }\, 
\overline{B^N(x,r)}\subset D\,\}.\]

From above it follows that $H\ne\emptyset$. It is  easy to see that $H$ is open. To show that $H$ is closed in $D$, take a sequence $x_j\in H$, 
$j=1,2,\dots $, $x_j\rightarrow x^*$ as $j\rightarrow +\infty $. Take $r^*>0$ such that $\overline{B^N(x^*,r^*)}\subset D$. Choose 
$x_{j_0}\in H\cap B^N(x^*,\frac{r^*}{2})$. Then there is $r_{j_0}>0$ such that $u(x)\in {\mathbb{R}}$ for almost every $x\in B^N(x_{j_0},r_{j_0})$. 
Hence also 
$u(x)\in {\mathbb{R}}$ for almost every $x\in B^N(x_{j_0},r_{j_0})\cap B^N(x^*,\frac{r^*}{2})$. Let $A$ be the set of such points. Clearly $A$ is of positive 
measure. Choose $\hat x\in A$. Then 
\[-\infty <u(\hat x)\leq \frac{K}{\nu _N\,\left(\frac{r^*}{2}\right)^N}\int\limits_{B^N(\hat x,\frac{r^*}{2})}u(y)\,dm_N(y)<+\infty .\]
Hence  $u(x)\in {\mathbb{R}}$ for almost every $x\in B^N(\hat x,\frac{r^*}{2})$. On the other hand, $x^*\in B^N(\hat x,\frac{r^*}{2})$. Thus $x^*$ has a 
neighborhood $B^N(x^*,\delta )\subset B^N(\hat x,\frac{r^*}{2})$ such that $u(x)\in {\mathbb{R}}$ for almost every $x\in B^N(x^*,\delta )$. 
 Hence $x^*\in H$.
Since $D$ is connected, $H=D$.

To show that $u\in {\mathcal{L}}^1_{\textrm{loc}}(D)$, take $\overline{B^N(x_0,r_0)}\subset D$ arbitrarily.  We know that for each $x\in \overline{B^N(x_0,r_0)}$ 
there is $r_x>0$ such that  $u(y)\in {\mathbb{R}}$ for almost every $y\in B^N(x,r_x)$.  Then
\[B^N(x,r_x),\quad x\in \overline{B^N(x_0,r_0)},\]
is an open cover of the compact set $\overline{B^N(x_0,r_0)}$. Hence we find a finite subcover
\[B^N(x_1,r_1), B^N(x_2,r_2), \dots , B^N(x_k,r_k).\]
We may suppose that $u(x_j)\in {\mathbb{R}}$ for  all $j=1,2,\dots ,k.$ One achieves this just replacing $x_j$, if necessary, by a nearby point $x^*_j$, and possibly 
increasing $r_j$ a little bit, $j=1,2,  \dots ,k$. Then 
\[-\infty <u(x_j)\leq \frac{K}{\nu _N\,r_j^N}\int\limits_{B^N(x_j,r_j)}u(y)dm_N(y), \quad j=1,2,\dots ,k.\]
By (iii) above, there is a constant $C>0$ such that $u(y)\leq u^+(y)\leq C<+\infty $ for each \mbox{$y\in \overline{B^N(x_1,r_1)}\cup \overline{B^N(x_2,r_2)}
\cup \cdots \cup  \overline{B^N(x_k,r_k)}$.}
Therefore  $u(y)-C\leq 0$ for each $y\in \overline{B^N(x_1,r_1)}\cup \overline{B^N(x_2,r_2)}\cup \cdots \cup  \overline{B^N(x_k,r_k)}$, and we get 
\begin{align*}-\infty &<\sum\limits_{j=1}^{k} \nu _N\, r_j^N\, u(x_j)\leq 
\sum\limits_{j=1}^{k}K\cdot\int\limits_{B^N(x_j,r_j)}u(y)dm_N(y)\\
&\leq K\cdot\sum\limits_{j=1}^{k}\int\limits_{B^N(x_j,r_j)}(u(y)-C)dm_N(y)+K\cdot C\cdot \sum\limits_{j=1}^{k}m_N(B^N(x_j,r_j))\\
&\leq K\cdot\int\limits_{\overline{B^N(x_0,r_0)}}(u(y)-C)dm_N(y)+K\cdot C\cdot \sum\limits_{j=1}^{k}m_N(B^N(x_j,r_j))\\
&\leq K\cdot\int\limits_{B^N(x_0,r_0)}u(y)dm_N(y)-K\cdot C\,m_N(B^N(x_0,r_0))+K\cdot C\cdot \sum\limits_{j=1}^{k}m_N(B^N(x_j,r_j))\\
&<+\infty. 
\end{align*}
Thus
\begin{displaymath}-\infty <\int\limits_{B^N(x_0,r_0)}u(y)\, dm_N(y)<+\infty ,\end{displaymath}
and  the claim follows. \hfill \qed

\vspace{1ex}
\noindent{\textbf{Remark.}} It is easy to see that (vii) does not anymore hold for quasi-nearly subharmonic functions. As a counterexample serves the function 
 $u:\,{\mathbb{R}}^2\rightarrow [-\infty ,+\infty )$,
\begin{displaymath} u(x,y):=\begin{cases}-\infty , & {\textrm{when }}\, y\leq 0,\\
1, & {\textrm{when }}\, y>0,\end{cases}\end{displaymath}
which is $2$-quasi-nearly subharmonic, but surely not quasi-nearly subharmonic n.s.

\vspace{4ex}

\noindent{{\textbf{3. \,\, Separately subharmonic functions}}} 

\vspace{3ex}

\noindent{\textbf{3.1.}} A counterpart and a generalization  to Theorem~A is the following  simple result: 

\vspace{1ex}

\noindent{\textbf{Proposition~3.}}  \emph{ Let $\Omega $ be a domain in ${\mathbb{R}}^{m+n}$, \mbox{$m,n\geq 2$, and let $K_1, K_2\geq 1$.} 
Let $u:\, \Omega \rightarrow 
[-\infty ,+\infty )$ be a Lebesgue measurable function such that} 
 \begin{itemize}
\item[(a)] \emph{for each $y\in {\mathbb{R}}^n$ the function} 
\[\Omega (y)\ni x\mapsto u(x,y)\in [-\infty ,+\infty )\]
\emph{is $K_1$-quasi-nearly subharmonic,}
 \item[(b)] \emph{for almost every $x\in {\mathbb{R}}^m$ the function} 
\[\Omega (x)\ni y\mapsto u(x,y)\in [-\infty ,+\infty )\]
\emph{is $K_2$-quasi-nearly subharmonic,}
\item[(c)] \emph{there exists a non-constant permissible function}   $\psi :\, [0,+\infty ) \rightarrow [0,+\infty )$ 
\emph{such that} $\psi \circ  u^+\in {\mathcal{L}}_{\textrm{loc}}^1(\Omega )$.
\end{itemize} 
 \emph{Then} $u$ 
\emph{is $\frac{4^{m+n}\nu _{m+n}K_1K_2}{\nu _m\,\nu _n}$-quasi-nearly subharmonic in $\Omega $}.

\vspace{1ex}

\noindent\emph{Proof.} We begin by showing that $u^+$ is locally bounded above in $\Omega $.  Since permissible functions are continuous, 
$\psi \circ u^+$ is Lebesgue measurable. Take  $(a,b)\in \Omega $ and $R>0$ arbitrarily such that $\overline{B^{m+n}((a,b),R)}\subset \Omega $. 
Take  $(x_0,y_0)\in B^m(a,\frac{R}{4})\times B^n(b,\frac{R}{4})$ 
arbitrarily. Using  assumptions (a) and (b), Proposition~2~(v) and Proposition~1~(iv), Fubini's Theorem and assumption~(c), one obtains:
\begin{align*}(\psi \circ u^+)(x_0,y_0)&\leq \frac{K_1}{\nu _m\,(\frac{R}{4})^{m}}\int\limits_{B^{m}(x_0,\frac{R}{4})}\, (\psi \circ u^+)(x,y_0)\,dm_{m}(x)\\
&\leq \frac{K_1}{\nu _m\,(\frac{R}{4})^{m}}\int\limits_{B^{m}(x_0,\frac{R}{4})}\big[ \frac{K_2}{\nu _n\,(\frac{R}{4})^n}\int\limits_{B^n(y_0,\frac{R}{4})}\,
(\psi \circ u^+)(x,y)\,dm_n(y)\big]dm_{m}(x)\\
&\leq \frac{4^{m+n}\nu _{m+n}K_1K_2}{\nu _m\,\nu _n}\frac{1}{\nu _{m+n}\,R^{m+n}}\int\limits_{B^{m+n}((a,b),R)}\, (\psi \circ u^+)(x,y)\,dm_{m+n}(x,y)\\
&<+\infty .
\end{align*}
Thus  $\psi \circ u^+$ is  locally 
bounded above in $\Omega $. Using  then  properties of  permissible functions,  one sees easily that also $u^+$ is locally bounded above in $\Omega $,
 thus locally integrable in $\Omega $. 

Proceeding then as above, but now
$\psi$ replaced with the identity mapping, and $u^+$ replaced with $u_M$, for an arbitrary  $M\geq 0$, and choosing $(x_0,y_0)=(a,b)$, one sees that
\begin{displaymath}u_M(a,b)\leq \frac{4^{m+n}\nu_{m+n}\,K_1K_2}{\nu_m\, \nu_n}\cdot \frac{1}{\nu _{m+n}R^{m+n}}
\int\limits_{B^{m+n}((a,b),R)}u_M(x,y)\,dm_{m+n}(x,y).\end{displaymath} 
Thus $u$ is  $\frac{4^{m+n}\nu _{m+n}K_1K_2}{\nu _m \nu_n}$-quasi-nearly 
 subharmonic in $\Omega $.\hfill\qed

\vspace{1ex}

Next we give:

\vspace{1ex}

\noindent{\textbf{Theorem~1.}}  \emph{Let $\Omega $ be a domain in ${\mathbb{R}}^{m+n}$, \mbox{$m,n\geq 2$, and let $K\geq 1$.} 
Let $u:\, \Omega \rightarrow 
[-\infty ,+\infty )$ be a Lebesgue measurable function such that} 
 \begin{itemize}
\item[(a)] \emph{for each $y\in {\mathbb{R}}^n$ the function} 
\[\Omega (y)\ni x\mapsto u(x,y)\in [-\infty ,+\infty )\]
\emph{is $1$-quasi-nearly subharmonic,}
 \item[(b)] \emph{for almost every $x\in {\mathbb{R}}^m$ the function} 
\[\Omega (x)\ni y\mapsto u(x,y)\in [-\infty ,+\infty )\]
\emph{is $K$-quasi-nearly subharmonic n.s.,}
\item[(c)] \emph{there exists a non-constant permissible function}   $\psi :\, [0,+\infty ) \rightarrow [0,+\infty )$ \emph{such that} 
$\psi \circ  u^+\in {\mathcal{L}}_{\textrm{loc}}^1(\Omega )$.
\end{itemize} 
 \emph{Then} $u$ 
\emph{is $K$-quasi-nearly subharmonic n.s. in $\Omega $}.

\vspace{1ex}

\noindent\emph{Proof.}  By Proposition~3 u is quasi-nearly subharmonic in $\Omega $. Thus  $u^+$ is locally integrable in $\Omega $. 

It remains to show that for all  $(a,b)\in \Omega $ and $R>0$  such that $\overline{B^{m+n}((a,b),R)}\subset \Omega $,
\begin{displaymath}u(a,b)\leq \frac{K}{\nu _{m+n}R^{m+n}}\int\limits_{B^{m+n}((a,b),R)}u(x,y)dm_{m+n}(x,y).\end{displaymath}
To see this, we proceed in the following standard way, see e.g. [Her71, proof of Theorem~a), pp.~32-33]:
\begin{align*}&\frac{K}{\nu _{m+n}R^{m+n}}\int\limits_{B^{m+n}((a,b),R)}u(x,y)dm_{m+n}(x,y)\\
=&\frac{\nu _n}{\nu _{m+n}R^{m+n}}\int\limits_{B^m(a,R)}[(R^2-\mid x-a\mid ^2)^{\frac{n}{2}}\frac{K}{\nu _n(R^2-\mid x-a\mid ^2)^{\frac{n}{2}}}
\int\limits_{B^n(b,\sqrt{R^2-\mid x-a\mid ^2})}u(x,y)dm_n(y)]dm_m(x) \\
\geq &\frac{\nu_n}{\nu _{m+n}R^{m+n}}\int\limits_{B^m(a,R)}(R^2-\mid x-a\mid ^2)^{\frac{n}{2}}u(x,b)dm_m(x)\geq u(a,b).\end{align*}
Above we have used, in addition to the fact that, for almost every $x\in {\mathbb{R}}^m$,  the functions $u(x,\cdot )$  are $K$-quasi-nearly 
subharmonic n.s., also the following lemma. (The proof of the Lemma, see [Her71, proof of Theorem~2~a), p.~15], works also in our slightly more 
general situation.) 

\vspace{1ex}

\noindent{\textbf{Lemma.}} ([Her71, Theorem~2~a), p.~15])  \emph{Let $v$ be  nearly subharmonic (in the generalized sense, defined above) in 
a domain $U$ of ${\mathbb{R}}^{N}$, $N\geq 2$, $\psi \in {\mathcal{L}}^{\infty }({\mathbb{R}}^N)$, $\psi \geq 0$, $\psi (x)=0$ when $\mid x\mid \geq \alpha $ and $\psi (x)$ depends only on $\mid x\mid $.
 Then $\psi \star v\geq v$ and $\psi \star v$ is subharmonic in $U_{\alpha }$, provided $\int \psi (x)dm_N(x)=1$, where 
$U_{\alpha }=\{x\in U:\, \overline{B^N(x,\alpha )}\subset U\}$.}

\vspace{1ex}

\par
\quad \hfill \qed

\vspace{1ex}

Choosing $K=1$ in Theorem~1 we get:

\vspace{1ex} 

\noindent{\textbf{Corollary~1.}}  \emph{ Let $\Omega $ be a domain in ${\mathbb{R}}^{m+n}$, \mbox{$m,n\geq 2$.} 
Let $u:\, \Omega \rightarrow 
[-\infty ,+\infty )$ be a Lebesgue measurable function such that} 
 \begin{itemize}
\item[(a)] \emph{for each $y\in {\mathbb{R}}^n$ the function} 
\[\Omega (y)\ni x\mapsto u(x,y)\in [-\infty ,+\infty )\]
\emph{is nearly subharmonic,}
 \item[(b)] \emph{for almost every $x\in {\mathbb{R}}^m$ the function} 
\[\Omega (x)\ni y\mapsto u(x,y)\in [-\infty ,+\infty )\]
\emph{is nearly subharmonic,}
\item[(c)] \emph{there exists a non-constant permissible function}   $\psi :\, [0,+\infty ) \rightarrow [0,+\infty )$ \emph{such that} 
$\psi \circ  u^+\in {\mathcal{L}}_{\textrm{loc}}^1(\Omega )$.
\end{itemize} 
 \emph{Then} $u$ 
\emph{is nearly subharmonic in $\Omega $}.

\vspace{1ex}

Before giving two further corollaries to Theorem~1, we state a more or less well-known measurability result. Our proof is an improved version of 
the proof of [Rii84, Lemma~3.2, pp. 103-104]. For a related result, with a different proof, see [Ar66, Lemma~1, p.~624].

\vspace{1ex}

\noindent{\textbf{Lemma.}}  \emph{ Let $U$ be a domain in ${\mathbb{R}}^{m}$ and $V$ be a domain in ${\mathbb{R}}^{n}$ , \mbox{$m,n\geq 1$.} 
Let $v:\, U\times V \rightarrow 
[-\infty ,+\infty )$ be such that} 
 \begin{itemize}
\item[(a)] \emph{for each $y\in V$ the function} 
\[U\ni x\mapsto v(x,y)\in [-\infty ,+\infty )\]
\emph{is  Lebesgue integrable and, for almost every $y\in V$ and every $x\in  U$,}
\[\frac{1}{\nu _mr^m}\int\limits_{B^m(x,r)}v(z,y)dm_m(z)\rightarrow v(x,y)\]
\emph{as $r\rightarrow 0$,}
 \item[(b)] \emph{for each $x\in U$ the function} 
\[V\ni y\mapsto v(x,y)\in [-\infty ,+\infty )\]
\emph{is upper semicontinuous.}
\end{itemize} 
 \emph{Then $v$ is Lebesgue measurable.} 

\vspace{1ex}

\noindent\emph{Proof.} It is clearly sufficient to show, that for any $M\geq 0$, the function $v_M=\max\{v,-M\}+M$ is measurable. Thus we may suppose that 
$v\geq 0$.

For each $r>0$ define $h_r:\, U\times V\rightarrow {\mathbb{R}}$,
\begin{displaymath}h_r(x,y)=\frac{1}{\nu _mr^m}\int\limits_{B^m(x,r)}v(z,y)dm_m(z).\end{displaymath}
Here we understand that $v(z,y)=0$ whenever $(x,y)\in  ({\mathbb{R}}^m\times {\mathbb{R}}^n)\setminus (U\times V)$.

To see that for each $x\in U$ the function $h_r(x,\cdot )$ is measurable, we proceed as follows. Write $v_k=\min\{v,k\}$ and 
$h_r^k:\, U\times V\rightarrow {\mathbb{R}}$,
\begin{displaymath}h^k_r(x,y)=\frac{1}{\nu _mr^m}\int\limits_{B^m(x,r)}v_k(z,y)dm_m(z),\end{displaymath}
$k=1,2,\dots $. Take $y\in  V$ and a sequence $y_j\rightarrow y$, $y_j\in  V$, $j=1,2,\dots$, arbitrarily. Using Fatou's Lemma one gets,
\begin{align*}\limsup_{j\rightarrow +\infty }h^k_r(x,y_j)=&\limsup_{j\rightarrow +\infty }\frac{1}{\nu _mr^m}\int\limits_{B^m(x,r)}v_k(z,y_j)dm_m(z)\\
 \leq &\frac{1}{\nu _mr^m}\int\limits_{B^m(x,r)}\limsup_{j\rightarrow +\infty }v_k(z,y_j)dm_m(z)\\
 \leq &\frac{1}{\nu _mr^m}\int\limits_{B^m(x,r)}v_k(z,y)dm_m(z)=h^k_r(x,y).
\end{align*}
Thus $h^k_r(x,\cdot )$ is upper semicontinuous in $V$. Using then Lebesgue Monotone Convergence Theorem one sees that, for each $x\in U$,
\begin{align*}\lim_{k\rightarrow +\infty }h^k_r(x,y)=&\lim_{k\rightarrow +\infty }\frac{1}{\nu _mr^m}\int\limits_{B^m(x,r)}v_k(z,y)dm_m(z)\\
 = &\frac{1}{\nu _mr^m}\int\limits_{B^m(x,r)}\lim_{k\rightarrow +\infty }v_k(z,y)dm_m(z)\\
 = &\frac{1}{\nu _mr^m}\int\limits_{B^m(x,r)}v(z,y)dm_m(z)=h_r(x,y).
\end{align*}
Hence $h_r(x,\cdot )$ is measurable.

To see that for each $y\in V$ the function $h_r(\cdot ,y)$ is continuous in $U$, we proceed as follows. Take $x,x_0\in  U$ arbitrarily. Then
\begin{align*} \mid h_r(x,y)-h_r(x_0,y)\mid =&\mid \frac{1}{\nu _mr^m}\int\limits_{B^m(x,r)}v(z,y)dm_m(z)-
\frac{1}{\nu _mr^m}\int\limits_{B^m(x_0,r)}v(z,y)dm_m(z)\mid \\
\leq &\int\limits_{B^m(x,r)\triangle B^m(x_0,r)}\mid v(z,y)\mid dm_m(z),
\end{align*}
and the continuity follows. 

From a  classical result, originally due to Lebesgue, it follows then that $h_r$ is measurable. Since for almost every $y\in V$ and every $x\in U$,
\begin{displaymath}\lim_{r\rightarrow 0}h_r(x,y)=v(x,y),\end{displaymath}
we see that $v$ is measurable, concluding the proof.\hfill \qed

\vspace{1ex}

Next a continuity result: 

\vspace{1ex}

\noindent{\textbf{Corollary~2.}}  \emph{ Let $\Omega $ be a domain in ${\mathbb{R}}^{m+n}$, \mbox{$m,n\geq 2$, and let $K\geq 1$.} 
Let $u:\, \Omega \rightarrow 
[-\infty ,+\infty )$ be such that} 
 \begin{itemize}
\item[(a)] \emph{for each $y\in {\mathbb{R}}^n$ the function} 
\[\Omega (y)\ni x\mapsto u(x,y)\in [-\infty ,+\infty )\]
\emph{is  nearly subharmonic, and, for almost every $y\in {\mathbb{R}}^n$, subharmonic,}
 \item[(b)] \emph{for each $x\in {\mathbb{R}}^m$ the function} 
\[\Omega (x)\ni y\mapsto u(x,y)\in [-\infty ,+\infty )\]
\emph{is upper semicontinuous, and, for almost every $x\in {\mathbb{R}}^m$, $K$-quasi-nearly subharmonic n.s.,}
\item[(c)] \emph{there exists a non-constant permissible function}   $\psi :\, [0,+\infty ) \rightarrow [0,+\infty )$ \emph{such that} 
$\psi \circ  u^+\in {\mathcal{L}}_{\textrm{loc}}^1(\Omega )$.
\end{itemize} 
 \emph{Then for every $(a,b)\in \Omega $},
\begin{displaymath}\limsup_{(x,y)\rightarrow (a,b)}u(x,y)\leq K\,u^+(a,b).\end{displaymath} 

\vspace{1ex}
\noindent\emph{Proof.} By the above Lemma $u$ is measurable. By Theorem~1  above, $u$ and thus also $u^+$ are $K$-quasi-nearly subharmonic n.s. Thus $u^+$ is locally bounded above in $\Omega $. 
Since also $u^+$ satisfies the assumptions 
of the corollary, it is sufficient to show that for any 
$(a,b)\in \Omega $, 
\begin{displaymath}\limsup_{(x,y)\rightarrow (a,b)}u^+(x,y)\leq K\, u^+(a,b).\end{displaymath}
Take $(a,b)\in \Omega $ and $R_1>0$ and $R_2>0$ arbitrarily such that $\overline{B^m(a,R_1)\times B^n(b,R_2)}\subset \Omega $. Choose an arbitrary $\lambda \in {\mathbb{R}}$ such that $u^+(a,b)<\lambda $. Since $u^+(a,\cdot )$ is upper 
semicontinuous, we find $R_2'$, $0<R_2'<R_2$, such that 
\begin{displaymath}\frac{1}{\nu_nR_2'^n}\int\limits_{B^n(b,R_2')}u^+(a,y)dm_n(y)<\lambda .\end{displaymath}
Using the fact that, for almost every $y\in {\mathbb{R}}^n$, the function $u^+(\cdot ,y)$, is subharmonic, we get
\begin{displaymath}\frac{1}{\nu_mr^m}\int\limits_{B^m(a,r)}u^+(x,y)dm_m(x)\rightarrow u^+(a,y) \, {\textrm{ as }}\, r\rightarrow 0.
\end{displaymath}
Since $u^+$ is   locally bounded above, one can use Lebesgue Monotone Convergence Theorem. Thus  we find
$R'_1$, $0<R'_1<R_1$, such that  
\begin{displaymath}\frac{1}{\nu _nR_2'^n}\int\limits_{B^n(b,R_2')}[\frac{1}{\nu_mR_1'^m}\int\limits_{B^m(a,R_1')}u^+(x,y)dm_m(x)]dm_n(y)<\lambda .
\end{displaymath}
Choose $r_1$, $0<r_1<R_1'$, and $r_2$, $0<r_2<R_2'$, arbitrarily. Then for each $(x,y)\in B^m(a,r_1)\times B^n(b,r_2)$,
\begin{align*}u^+(x,y)&\leq \frac{1}{\nu _m(R_1'-r_1)^m}\int\limits_{B^m(x,R_1'-r_1)}u^+(\xi ,y)dm_m(\xi )\\
                      &\leq \frac{1}{\nu _m(R_1'-r_1)^m}\int\limits_{B^m(x,R_1'-r_1)}[\frac{K}{\nu _n(R_2'-r_2)^n}\int\limits_{B^n(y,R_2'-r_2)}
u^+(\xi ,\eta )dm_n(\eta )] dm_m(\xi )\\
&\leq \frac{K}{\nu _n(R_2'-r_2)^n}\int\limits_{B^n(y,R_2'-r_2)}[\frac{1}{\nu _m(R_1'-r_1)^m}\int\limits_{B^m(x,R_1'-r_1)}
u^+(\xi ,\eta )dm_m(\xi )] dm_n(\eta )\\
&\leq \left( \frac{R_1'}{R_1'-r_1}\right)^m\cdot \left(\frac{R_2'}{R_2'-r_2}\right)^n \cdot \frac{K}{\nu _nR_2'^n}\int\limits_{B^n(b,R_2')}
[\frac{1}{\nu _mR_1'^m}\int\limits_{B^m(a,R_1')}u^+(\xi ,\eta )dm_m(\xi )]dm_n(\eta )\\
&\leq \left(\frac{R_1'}{R_1'-r_1}\right)^m\cdot \left(\frac{R_2'}{R_2'-r_2}\right)^n \cdot K \cdot \lambda .
\end{align*}
Sending then $r_1\rightarrow 0$, $r_2\rightarrow 0$, one gets 
\begin{displaymath}\limsup_{(x,y)\rightarrow (a,b)}u^+(x,y)\leq K\cdot \lambda ,\end{displaymath}
concluding the proof. \hfill \qed 

\vspace{1ex}

\noindent{\textbf{Remark.}} The above proof is essentially the same as [Ri89, part of the  proof of Theorem~1, pp.~70-71], where we gave a short proof for the upper 
semicontinuity of a locally integrable separately subharmonic function.

Our last corollary  improves our Theorem~A, and thus also the cited results of Lelong, Avanissian and Arsove.

\vspace{1ex} 

\noindent{\textbf{Corollary~3.}}  \emph{ Let $\Omega $ be a domain in ${\mathbb{R}}^{m+n}$, \mbox{$m,n\geq 2$.} 
Let $u:\, \Omega \rightarrow 
[-\infty ,+\infty )$ be  such that} 
 \begin{itemize}
\item[(a)] \emph{for each  $y\in {\mathbb{R}}^n$ the function} 
\[\Omega (y)\ni x\mapsto u(x,y)\in [-\infty ,+\infty )\]
\emph{is nearly subharmonic, and, for almost every $y\in {\mathbb{R}}^n$, subharmonic,}
 \item[(b)] \emph{for each $x\in {\mathbb{R}}^m$ the function} 
\[\Omega (x)\ni y\mapsto u(x,y)\in [-\infty ,+\infty )\]
\emph{is upper semicontinuous, and, for almost every $x\in {\mathbb{R}}^m$,  subharmonic,}
\item[(c)] \emph{there exists a non-constant permissible function}   $\psi :\, [0,+\infty ) \rightarrow [0,+\infty )$ \emph{such that} 
$\psi \circ  u^+\in {\mathcal{L}}_{\textrm{loc}}^1(\Omega )$.
\end{itemize} 
 \emph{Then} $u$ 
\emph{is  subharmonic in $\Omega $}.

\vspace{1ex}

\noindent\emph{Proof.} By the above Lemma, $u$ is measurable. By Corollary~1, $u$ and thus also $u_M=\max\{u,-M\}+M$, for each $M\geq 0$,    are nearly subharmonic.  
It is clearly sufficient to show that $u_M$ is upper semicontinuous. Since $u_M$ satisfies the assumptions of  Corollary~2, the claim 
follows. \hfill\qed

\vspace{1ex} 

\noindent{\textbf{Remark.}} Compare the above short proof for Corollary~3 to  the  proofs of the previous results 
  [Le45, Théorème~1~bis, pp.~308-315],  [Av61, proofs of Proposition~10 and Théorème~9, pp.~137-140], see also  [Her71, proof of Theorem, pp.~31-32],
and [Ar66, proof of Theorem~1, pp.~624-625]. It is, however, worth while to point out that still another, 
a new and perhaps an even more direct proof for Corollary~3  exists:   

\vspace{1ex}

\noindent\emph{A direct proof for Corollary~3}. Since by Corollary~1 $u$ is nearly subharmonic, it remains to show that $u$ is 
upper semicontinuous. But this follows at once from the following lemma (whose proof, see [Her71, pp.~34-35], uses the simple \emph{iterated}  mean  value inequality, 
but \emph{not} the general mean value inequality!):

\vspace{1ex}

\noindent{\textbf{Lemma.}} ([Her71, Proposition~2, p.~34]) \emph{ Let $\Omega $ be a domain in ${\mathbb{R}}^{m+n}$, \mbox{$m,n\geq 2$.} 
Let $u:\, \Omega \rightarrow 
[-\infty ,+\infty )$ be nearly subharmonic and  such that} 
 \begin{itemize}
\item[(a)] \emph{for almost every  $y\in {\mathbb{R}}^n$ the function} 
\[\Omega (y)\ni x\mapsto u(x,y)\in [-\infty ,+\infty )\]
\emph{is subharmonic,}
 \item[(b)] \emph{for almost every $x\in {\mathbb{R}}^m$ the function} 
\[\Omega (x)\ni y\mapsto u(x,y)\in [-\infty ,+\infty )\]
\emph{is nearly subharmonic.}
\end{itemize} 
 \emph{Then the same subharmonic function} $u^*$ \emph{is obtained by regularization either with respect to $(x,y)$ or with respect to $y$ only 
(owing to dissymmetry in the assumptions).}
\par
\hfill\qed

\vspace{1ex}

\noindent{\textbf{3.2. Remark.}} Unlike in Theorem~A, the  measurability assumptions in Proposition~3 and Theorem~1 are now necessary. 
With the aid of Sierpinski's  nonmeasurable function, given e.g. in  [Ru79, 7.9 (c), pp.~152-153], one easily constructs 
a nonmeasurable, separately quasi-nearly subharmonic  function $u:\,{\mathbb{C}}^2\rightarrow [1,2]$. Indeed, let $\tilde{Q}=\{0\}\times Q\times \{0\}
\subset {\mathbb{R}}\times {\mathbb{R}}^2\times {\mathbb{R}}$, where $Q\subset {\mathbb{R}}^2$ is the set of Sierpinski, see [Ru79, pp.~152-153, and 
Theorem~7.2, p.~146]. Then 
the function $v(z_1,z_2)=v(x_1,y_1,x_2,y_2):=2-\chi _{\tilde{Q}}(z_1,z_2)$ is clearly nonmeasurable, but still separately quasi-nearly  subharmonic, 
see Example~1 in {\textbf{2.3}} above. 

 \vspace{4ex}
\noindent{{\textbf{4. \,\, The result of Arsove}}} 

\vspace{3ex}

\noindent{\textbf{4.1.}} Next we give a slight generalization to a result of Arsove and of  Cegrell and Sadullaev, Theorem~B above. For the 
short proof given below, compare  also the proofs of the special case results [Ar66, proof of Theorem~2, p.~625] and [Ri06$_2$, proof of Theorem~B, 4.2]).  

\vspace{1ex}

\noindent{\textbf{Theorem~2.}} \emph{Let $\Omega $ be a domain in ${\mathbb{R}}^{m+n}$, \mbox{$m,n\geq 2$, and $K\geq 1$.} 
Let $u:\, \Omega \rightarrow 
{\mathbb{R}}$ be  such that} 
 \begin{itemize}
\item[(a)] \emph{for each $y\in {\mathbb{R}}^n$ the function} 
\[\Omega (y)\ni x\mapsto u(x,y)\in {\mathbb{R}}\]
\emph{is K-quasi-nearly subharmonic n.s.,}
 \item[(b)] \emph{for each $x\in {\mathbb{R}}^m$ the function} 
\[\Omega (x)\ni y\mapsto u(x,y)\in {\mathbb{R}}\]
\emph{is harmonic,}
\item[(c)] \emph{there is a nonnegative function } $\varphi \in {\mathcal{L}}^1_{\textrm{loc}}(\Omega )$ \emph{such that} $-\varphi \leq u$.
\end{itemize}
\emph{Then} $u$ 
\emph{is K-quasi-nearly subharmonic n.s. in $\Omega $}.

\vspace{1ex}

\noindent\emph{Proof.} It is easy to see that $u$ is Lebesgue measurable. (See the end of the proof of the  above Lemma in {\textbf{3.1}}.) By Theorem~1 above it is sufficient 
to show that $u^+\in {\mathcal{L}}^1_{\textrm{loc}}(\Omega )$. Write $v:=u+\varphi $. Then $0\leq u^+\leq v$.  It is sufficient to show that
$v\in {\mathcal{L}}^1_{\textrm{loc}}(\Omega )$. Let $(a,b)\in \Omega $ and $R>0$ be 
arbitrary such that $\overline{B^m(a,R)\times B^n(b,R)}\subset \Omega $. 
Then
\begin{displaymath}\begin{split}0&\leq \frac{K}{m_{m+n}(B^m(a,R)\times B^n(b,R))}\int\limits_{B^m(a,R)\times B^n(b,R)}v(x,y) dm_{m+n}(x,y)\\
  =&\frac{K}{\nu _m\, R^m}\int\limits_{B^m(a,R)}\{\frac{1}{\nu _n\,R^n}\int\limits_{B^n(b,R)} [u(x,y)+\varphi (x,y)] dm_n(y)\}dm_m(x)\\
=&\frac{K}{\nu _m\, R^m}\int\limits_{B^m(a,R)}[\frac{1}{\nu _n\,R^n}\int\limits_{B^n(b,R)} u(x,y)dm_n(y)+
\frac{1}{\nu _n\,R^n}\int\limits_{B^n(b,R)} \varphi (x,y) dm_n(y)]dm_m(x)\\
=&\frac{K}{\nu _m\, R^m}\int\limits_{B^m(a,R)}[u(x,b)+
\frac{1}{\nu _n\,R^n}\int\limits_{B^n(b,R)} \varphi (x,y) dm_n(y)]dm_m(x)\\
=&\frac{K}{\nu _m\, R^m}\int\limits_{B^m(a,R)}u(x,b)dm_m(x)+\frac{K}{\nu _mR^m}\int\limits_{B^m(a,R)}[
\frac{1}{\nu _n\,R^n}\int\limits_{B^n(b,R)} \varphi (x,y) dm_n(y)]dm_m(x)\\
=&\frac{K}{\nu _m\, R^m}\int\limits_{B^m(a,R)}u(x,b)dm_m(x)+\frac{K}{m_{m+n}(B^m(a,R)\times B^n(b,R))}\int\limits_{B^m(a,R)\times B^n(b,R)}
\varphi (x,y) dm_{m+n}(x,y)\\
<&+\infty .
\end{split}\end{displaymath}
\hfill \qed

\vspace{1ex}

The following corollary improves the result of Arsove and Cegrell and Sadullaev, Theorem~B above:

\vspace{1ex}

\noindent{\textbf{Corollary.}} \emph{Let $\Omega $ be a domain in ${\mathbb{R}}^{m+n}$, \mbox{$m,n\geq 2$.} 
Let $u:\, \Omega \rightarrow 
{\mathbb{R}}$ be  such that} 
 \begin{itemize}
\item[(a)] \emph{for each $y\in {\mathbb{R}}^n$ the function} 
\[\Omega (y)\ni x\mapsto u(x,y)\in {\mathbb{R}}\]
\emph{is nearly subharmonic, and, for almost every $y\in {\mathbb{R}}^n$, subharmonic,}
 \item[(b)] \emph{for each $x\in {\mathbb{R}}^m$ the function} 
\[\Omega (x)\ni y\mapsto u(x,y)\in {\mathbb{R}}\]
\emph{is  harmonic,}
\item[(c)] \emph{there is a nonnegative function } $\varphi \in {\mathcal{L}}^1_{\textrm{loc}}(\Omega )$ \emph{such that} $-\varphi \leq u$.
\end{itemize}
\emph{Then} $u$ 
\emph{is subharmonic in $\Omega $}.

\vspace{1ex}

\noindent\emph{Proof.} By Theorem~2, $u\in {\mathcal{L}^1_{\textrm{loc}}}(\Omega )$ and thus also $u^+\in {\mathcal{L}^1_{\textrm{loc}}}(\Omega )$. Thus the claim follows from the above Corollary~3.\hfill \qed

\vspace{2ex}

\noindent{\textbf{4.2.}}  Also the following result gives a partial generalization to Theorem~B above. It   generalizes slightly  
[CS93, Corollary, p.~82], too. 
 
\vspace{1ex}

\noindent{\textbf{Theorem~3.}}  \emph{ Let $\Omega $ be a domain in ${\mathbb{R}}^{m+n}$, \mbox{$m,n\geq 2$.} 
Let $u:\, \Omega \rightarrow [0,+\infty )$ be such that} 
 \begin{itemize}
\item[(a)] \emph{for each $y\in {\mathbb{R}}^n$ the function} 
\[\Omega (y)\ni x\mapsto u(x,y)\in [0,+\infty )\]
\emph{is quasi-nearly subharmonic,}
 \item[(b)] \emph{for each $x\in {\mathbb{R}}^m$ the function} 
\[\Omega (x)\ni y\mapsto u(x,y)\in [0,+\infty )\]
\emph{is a $\lambda $-Harnack function, where $\lambda \in (0,1)$ is fixed.}
\end{itemize}
\emph{Then} $u$ 
\emph{is quasi-nearly subharmonic in $\Omega$}.

\vspace{1ex}

\noindent{\emph{Proof}}. It is well-known that $u$ is Lebesgue measurable, see the  proof of the previous Lemma, above in {\textbf{3.1.}} Let $(a,b)\in \Omega $ and $R>0$ be arbitrary such that $\overline{B^{m+n}((a,b),R)}\subset \Omega $. Choose 
$(x_0,y_0)\in B^m(a,\frac{R}{4})\times B^n(b,\frac{R}{4})$ arbitrarily. Since $u(\cdot ,y_0)$ is quasi-nearly subharmonic, one has, for some $K\geq 1$,
\begin{equation*} u(x_0,y_0)\leq \frac{K}{\nu _m(\frac{R}{4})^m}\int\limits_{B^m(x_0,\frac{R}{4})}u(x,y_0)\, dm_m(x).\end{equation*}
On the other hand, since the functions $u(x,\cdot )$, $x\in B^m(a,\frac{R}{2})$, are Harnack functions in $B^n(b,\frac{R}{2})$, there is a constant
$C=C(n,\lambda ,C_\lambda ,R)$ (here $\lambda $ and $C_\lambda $ are the constants in {\textbf{2.2.}}) such that 
\[\frac{1}{C}\leq \frac{u(x,y_0)}{u(x,b)}\leq C\]
for each $x\in B^m(a,\frac{R}{2})$. See e.g. [ABR01, proof of 3.6, pp.~48--49]. Therefore
 \begin{align*} u(x_0,y_0)&\leq \frac{K}{\nu _m(\frac{R}{4})^m}\int\limits_{B^m(x_0,\frac{R}{4})}\, C \, u(x,b)\, dm_m(x)\\
&\leq \frac{C\cdot K}{\nu _m(\frac{R}{4})^m}\int\limits_{B^m(a,\frac{R}{2})}\,  u(x,b)\, dm_m(x)\\
&\leq \frac{4^mC\cdot K}{\nu _mR^m}\int\limits_{B^m(a,\frac{R}{2})}\, u(x,b)\, dm_m(x)<+\infty .\end{align*}
Thus $u$ is locally bounded above in $B^m(a,\frac{R}{4})\times B^n(b,\frac{R}{4})$, and therefore the result follows from Proposition~1~(v) and 
Proposition~3 above.\hfill \qed

\vspace{4ex}
\noindent{{\textbf{5.\,\, The result of Ko\l odziej and Thorbi\"ornson}}}

\vspace{3ex}

\noindent\textbf{5.1.} In our  generalization to the cited result of Ko\l odziej and Thorbi$\ddot {\textrm{o}}$rnson, we use the generalized 
Laplacian, defined with the aid of the Blaschke-Privalov operators, see e.g. [Sz33], [Sa41], [Ru50], [Br69], [Pl70], [Sh71] and [Sh78]. Let  $D$ be a domain in 
${\mathbb{R}}^{N}$,
 \mbox{$N\geq 2$,} and  $f:\, D \rightarrow {\mathbb{R}}$, $f\in {\mathcal{L}}^1_{\textrm{loc}}(D)$. We write 
\begin{displaymath}\begin{split}
\Delta_*f(x):&=\liminf_{r\rightarrow 0}\frac{2(N+2)}{r^2}\cdot\big[\frac{1}{\nu _Nr^N}\int\limits_{B^N(x,r)}f(x')dm_N(x')-f(x)\big],\\
\Delta^* f(x):&=\limsup_{r\rightarrow 0}\frac{2(N+2)}{r^2}\cdot\big[\frac{1}{\nu _Nr^N}\int\limits_{B^N(x,r)}f(x')
dm_N(x')-f(x)\big].
\end{split}\end{displaymath}
If $\Delta _*f(x)= \Delta ^*f(x),$ then  write $\Delta f(x):= \Delta _*f(x)=\Delta ^*f(x)$. 

If $f\in {\mathcal{C}}^2(D)$,
then 
\begin{displaymath}\Delta f(x)=(\sum\limits_{j=1}^{N}\frac{\partial^2f}{\partial x_j^2})(x),
\end{displaymath}
the standard Laplacian with respect to the variable $x=(x_1,x_2,\dots ,x_N)$. 
More generally, if $x\in D$ and 
$f\in t^1_2(x)$, i.e. $f$ has an ${\mathcal{L}}^1$ total differential  of order $2$ at $x$, 
then $\Delta f(x)$ equals with the pointwise Laplacian of $f$ at $x$, i.e.
\begin{displaymath}\Delta f(x)=\sum\limits_{j=1}^{N}D_{jj}f(x).\end{displaymath}
Here $D_{jj}f$ represents a generalization to the usual $\frac{\partial^2f}{\partial x_j^2}$,  $j=1,2, \dots ,N$.  See e.g. [CZ61, p.~172], [Sh56, p.~498], 
[Sh71, p.~369] and [Sh78, p.~29].

Recall that there are functions which are not ${\mathcal{C}}^2$ but 
for which the generalized Laplacian is nevertheless continuous:

\vspace{1ex}

\noindent{\textbf{Example~1.}} ([Sh78, p.~31]) The function $f:\,{\mathbb{R}}^N\rightarrow {\mathbb{R}}$,
\begin{displaymath} f(x):=\begin{cases}-1, & {\textrm{when }}\, x_N<0,\\
0, & {\textrm{when }}\, x_N=0,\\
1, & {\textrm{when }}\, x_N> 0, \end{cases}\end{displaymath}
is non-continuous, but nevertheless  $\Delta f(x)=0$ for all $x\in {\mathbb{R}}^N$.

\vspace{1ex}

\noindent{\textbf{Example~2.}} Let $1\leq k\leq N$ and let $E=\{(0,0,\dots ,0)\}\times {\mathbb{R}}^{N-k}\subset {\mathbb{R}}^N$. Let $0<\lambda \leq 1$.
Define $f: \, {\mathbb{R}}^N\rightarrow {\mathbb{R}}$, 
\[f(x)=f(x_1,x_2,\dots ,x_k,x_{k+1},\dots ,x_N):=\left(\sqrt{x_1^2+x_2^2+\cdots +x_k^2}\right)^{\lambda }.\] 
Then $f$ is continuous and subharmonic in ${\mathbb{R}}^N$, but not in ${\mathcal{C}}^1({\mathbb{R}}^N)$. Nevertheless,  $\Delta f$  is defined everywhere
in ${\mathbb{R}}^N$,  equals $+\infty $ 
in $E$, and  continuous in ${\mathbb{R}}^N$, in $E$  in the extended sense, with respect to the spherical metric:
\[q(a,b):=\begin{cases}\frac{\mid a-b\mid }{\sqrt{1+a^2}\,\sqrt{1+b^2}} &{\textrm{ when }}a, b\in [0,+\infty ),\\
\frac{1}{\sqrt{1+a^2}} &{\textrm{ when }}a\in [0,+\infty ){\textrm{ and }}b=+\infty .\end{cases}\]
Observe that $([0,+\infty ],q)$ is a complete metric space. 

\vspace{1ex}

If $f$ is subharmonic on $D$, it follows from [Sa41, p.~451] (see also [Ru50, Lemma~2.2, p.~280]) that    
$\Delta f(x):=\Delta_{*}f(x)=\Delta^*f(x)\in {\mathbb{R}}$ 
for almost every $x\in D$. 

Below the following notation is used. Let  $\Omega $ is a domain in ${\mathbb{R}}^{m+n}$, \mbox{$m,n\geq 2$}, and 
$u:\, \Omega \rightarrow {\mathbb{R}}$. If  $y\in {\mathbb{R}}^n$ is such that the function 
\begin{displaymath} \Omega (y)\ni x\mapsto f(x):=u(x,y)\in {\mathbb{R}}\end{displaymath}
is in ${\mathcal{L}}^1_{\textrm{loc}}(\Omega (y))$, then we write $\Delta _{1*}u(x,y):=\Delta_*f(x)$, 
$\Delta^* _1u(x,y):=\Delta^*f(x)$, and $\Delta _1u(x,y):=\Delta f(x)$.

\vspace{2ex}

\noindent\textbf{5.2.}  Then a generalization to the cited result Theorem~C of  Ko\l odziej and Thorbi$\ddot {\textrm{o}}$rnson [KT96, Theorem~1, p.~463]. 
Our result improves our 
previous result [Ri07$_2$, Theorem~1, p.~438]: Our assumptions (d) and (e)  are now essentially milder than previously. Consequently,  the present proof 
includes new ingredients, too.

\vspace{1ex}

\noindent{\textbf{Theorem~4.}}  \emph{ Let $\Omega $ be a domain in ${\mathbb{R}}^{m+n}$, \mbox{$m,n\geq 2$.} 
Let $u:\, \Omega \rightarrow 
{\mathbb{R}}$ be such that for each $(x',y')\in \Omega $ there is  $(x_0,y_0)\in \Omega $ and $r_1>0$, $r_2>0$  such that 
$(x',y')\in B^m(x_0,r_1)\times B^n(y_0,r_2)\subset \overline{B^m(x_0,r_1)\times B^n(y_0,r_2)}\subset \Omega $ and  such that the following 
conditions are satisfied:} 
 \begin{itemize}
\item[(a)] \emph{For each $y\in \overline{B^n(y_0,r_2)}$ the function} 
\[\overline{B^m(x_0,r_1)}\ni x\mapsto u(x,y)\in {\mathbb{R}}\]
\emph{is continuous, and subharmonic in $B^m(x_0,r_1)$.}
 \item[(b)] \emph{For each $x\in \overline{B^m(x_0,r_1)}$ the function} 
\[\overline{B^n(y_0,r_2)}\ni y\mapsto u(x,y)\in {\mathbb{R}}\]
\emph{is continuous, and harmonic in $B^n(y_0,r_2)$.}
\item[(c)] \emph{For each $y\in B^n(y_0,r_2)$ one has $\Delta _{1*}u(x,y)<+\infty $ for each $x\in B^m(x_0,r_1)$, possibly with the exception of 
a polar set in $B^m(x_0,r_1)$.}
\item[(d)] \emph{There is a set $H\subset B^n(y_0,r_2)$, dense in $B^n(y_0,r_2)$, and an open  set $K_1\subset B^m(x_0,r_1)$, dense in $B^m(x_0,r_1)$,
such that for each $y\in H$, for almost every $x\in B^m(x_0,r_1)$ and for any sequence $x_j\rightarrow x$, $x_j\in K_1$, $j=1,2,\dots$, 
the sequence $\Delta_{1*}u(x_j,y)$ is convergent in $([0,+\infty ],q)$.}
\item[(e)] \emph{There is an open  set $K_2\subset B^m(x_0,r_1)$, dense in $B^m(x_0,r_1)$, such that for each $y\in B^n(y_0,r_2)$, 
 for almost every  $x\in B^m(x_0,r_1)$, and for any  sequence \mbox{$x_j\rightarrow x$}, $x_j\in K_2$, 
$j=1, 2, \dots$, there is a subsequence $x_{j_l}$ such that  
\[\Delta _{1*}u(x_{j_l},y)\rightarrow \Delta _{1*}u(x,y)\]
in $([0,+\infty ],q)$ as $l\rightarrow +\infty $.} 
\end{itemize}
\emph{Then} $u$ 
\emph{is subharmonic in $\Omega $}.

\vspace{1ex}

\noindent\emph{Proof}. Choose $r'_1$,  $r'_2$ such that $0<r'_1<r_1$, $0<r'_2<r_2$, and  such that $(x',y')\in B^m(x_0,r'_1)\times B^n(y_0,r'_2)$.
It is sufficient to show that 
$u\mid B^m(x_0,r'_1)\times B^n(y_0,r'_2)$ is subharmonic. For the sake of convenience of notation, we change the roles of $r_j$ and $r'_j$, $j=1,2$.
We divide the proof into several steps.

\vspace{1ex}

\noindent{\textbf{Step~1.}} \emph{Construction of an auxiliar set $G$.}

\vspace{1ex}

For each $k\in {\mathbb{N}}$ write
\[A_k:=\{\, x\in \overline{B^m(x_0,r_1)}\, :\, -k\leq u(x,y)\leq k \quad {\textrm{for each}}\quad y\in  \overline{B^n(y_0,r_2)}\,\}.\] 
Clearly $A_k$ is closed, and 
\[\overline{B^m(x_0,r_1)}=\bigcup _{k=1}^{+\infty }A_k.\]
Write
\begin{displaymath}G:=\bigcup _{k=1}^{+\infty }{\textrm{int}} A_k.\end{displaymath}
It follows from Baire's Theorem  that $G$ is dense in $B^m(x_0,r_1)$.

\vspace{1ex}

\noindent{\textbf{Step~2.}} \emph{The functions $\Delta_{1r}u(x,\cdot )$ (see the definition below)}, $x\in G$, 
$0<r<r_x:={\textrm{dist}}(x,\overline{B^m(x_0,r_1)}\setminus G)$, \emph{are 
nonnegative and harmonic in $B^n(y_0,r_2)$.}

\vspace{1ex}

For each $(x,y)\in B^m(x_0,r_1)\times B^n(y_0,r_2)$ and  each $r$, $0<r<{\textrm{dist}}(x,\partial B^m(x_0,r'_1))$ (observe that 
${\textrm{dist}}(x,\partial B^m(x_0,r'_1))>r'_1-r_1>0$), write 
\begin{displaymath}\begin{split}\Delta _{1r}u(x,y)&:=\frac{2(m+2)}{r^2}\cdot \big[\frac{1}{\nu _m\, r^m}\int\limits_{B^m(x,r)}u(x',y)\, dm_m(x')-u(x,y)\big]\\
&=\frac{2(m+2)}{r^2}\cdot \frac{1}{\nu _m\, r^m}\int\limits_{B^m(0,r)}\big[u(x+x',y)-u(x,y)\big]\, dm_m(x').\end{split}\end{displaymath}    
Since $u(\cdot ,y)$ is subharmonic, $\Delta _{1r}u(x,y)$ is defined and nonnegative. Suppose then that $x\in G$ and $0<r<r_x$.
Since $\overline{B^m(x,r)}\subset G$ and $A_k\subset A_{k+1}$ for all $k=1,2, \dots$, $\overline{B^m(x,r)}\subset {\textrm{int}}A_N$ for some 
$N\in {\mathbb{N}}$. Therefore
\begin{displaymath}-N\leq u(x',y)\leq N\quad {\textrm{for all}} \quad x'\in B^m(x,r) \quad {\textrm{and}} \quad y\in B^n(y_0,r_2),\end{displaymath}
and hence
\begin{equation}-2N\leq u(x+x',y)-u(x,y)\leq 2N\quad {\textrm{for all}} \quad x'\in B^m(0,r)\quad  {\textrm{and}} \quad y\in B^n(y_0,r_2).\end{equation}
To show that $\Delta_{1r}u(x,\cdot )$ is continuous, pick an arbitrary sequence $y_j\rightarrow \tilde{y}_0$, $y_j, \tilde{y}_0\in B^n(y_0,r_2)$, 
$j=1,2, \dots$. Using then (3), Lebesgue Dominated Convergence Theorem and the continuity of $u(x, \cdot )$, one sees easily that 
$\Delta _{1r}u(x,\cdot )$ is continuous.

To show that $\Delta _{1r}u(x,\cdot )$ satisfies the mean value equality, take $\tilde{y}_0\in B^n(y_0,r_2)$ and $\rho >0$ 
such that $\overline{B^n(\tilde{y}_0,\rho )}\subset B^n(y_0,r_2)$. Because of (3) we can use Fubini's Theorem. Thus     
\begin{displaymath}\begin{split}\frac{1}{\nu _n\rho ^n}&\int\limits_{B^n(\tilde{y}_0,\rho )}\Delta _{1r}u(x,y) dm_n(y)\\
&=\frac{1}{\nu _n\rho ^n}\int\limits_{B^n(\tilde{y}_0,\rho )}\{\frac{2(m+2)}{r^2}\cdot \frac{1}{\nu _m
 r^m}\int\limits_{B^m(0,r)}\big[u(x+x',y)-u(x,y)\big] dm_m(x')\}dm_n(y)\\
&=\frac{2(m+2)}{r^2}\cdot \frac{1}{\nu _m
 r^m}\int\limits_{B^m(0,r)}\{\frac{1}{\nu _n\rho ^n}\int\limits_{B^n(\tilde{y}_0,\rho )}\big[u(x+x',y)-u(x,y)\big]dm_n(y)\} dm_m(x')\\
&=\frac{2(m+2)}{r^2}\cdot \frac{1}{\nu _m\, r^m}\int\limits_{B^m(0,r)}\big[u(x+x',\tilde{y}_0)-u(x,\tilde{y}_0)\big]\, dm_m(x')\\
&=\Delta _{1r}u(x,\tilde{y}_0).\end{split}\end{displaymath}     

\vspace{1ex}

\noindent{\textbf{Step~3.}} \emph{The functions  $\Delta_{1}u(x,\cdot ): B^n(y_0,r_2)\rightarrow {\mathbb{R}}$, $x\in G_1$ (see below for the definition 
of $G_1$),  are defined, nonnegative and harmonic.} 

\vspace{1ex}

By [Ru50, Lemma~2.2, p.~280] (see also [Sa41, p.~451]), for each $y\in B^n(y_0,r_2)$ there is a set $A(y)\subset B^m(x_0,r_1)$ such that 
$m_m(A(y))=0$ and that  
\begin{displaymath}\Delta_{1*}u(x,y)=\Delta_1^*u(x,y)=\Delta_1u(x,y)\in {\mathbb{R}}\end{displaymath}
for all $x\in B^m(x_0,r_1)\setminus A(y)$. 
We may clearly suppose that $H$ is countable, $H=\{\, y_k, \, k=1,2,\dots \, \}$.  
Write 
\[A:=\bigcup _{k=1}^{+\infty }A(y_k).\]
Since $G$ is open and dense in $B^m(x_0,r_1)$ and $A$ is of Lebesgue measure zero, also $G_1:=G\setminus A$ is dense in $B^m(x_0,r_1)$. 
Take $x\in G_1$ and a sequence  $r_j\rightarrow 0$, $0<r_j<r_x$, $j=1,2,\dots$, arbitrarily. Since $x\in G_1$, we see e.g. by  
[Her71, Corollary~3, p.~6] (or [AG01, Lemma~1.5.6, p.~16]), that the family 
\begin{displaymath} \Delta_ {1r_j}u(x,\cdot ), \, j=1,2,\dots ,\end{displaymath}
of nonnegative and harmonic functions in $B^n(y_0,r_2)$ is either equicontinuous and locally  uniformly bounded, or else 
\[\sup_{j=1,2,\dots}\Delta_{1r_j}u(x_j,\cdot )\equiv +\infty .\]
On the other hand, since $x\in G_1$, we know that
\begin{displaymath}\Delta_{1r_j}u(x,y_k)\rightarrow \Delta_1u(x,y_k)\in {\mathbb{R}}\, {\textrm{ as }}\, j\rightarrow +\infty ,\, 
{\textrm{ for each }}\, y_k\in H, \,k=1,2,\dots .\end{displaymath}   
Therefore,  by [V\"a71, Theorem~20.3, p.~68] and by [Her71, c), b), pp. 2, 3] (or [AG01,  Theorem~1.5.8, p.~17]), the limit 
\begin{displaymath}\Delta_1u(x,\cdot )=\lim_{j\rightarrow +\infty }\Delta_{1r_j}u(x,\cdot )\end{displaymath}
exists and defines a harmonic function in $B^n(y_0,r_2)$. Since the limit is clearly independent of the considered sequence $r_j$, the claim follows. 

\vspace{1ex}

\noindent{\textbf{Step~4.}} \emph{The function  $\Delta_{1}u(\cdot ,\cdot )\mid G_2\times B^n(y_0,r_2)$ 
has a continuous  extension $\tilde{\Delta}_{1}u(\cdot ,\cdot ): B\times B^n(y_0,r_0)\rightarrow {\mathbb{R}}$ 
(see below for the definitions of $G_2$ and $B$).  Moreover, the functions 
$\tilde{\Delta}_1u(x,\cdot )$, $x\in B$, are nonnegative and harmonic  in $B^n(y_0,r_2)$.}

\vspace{1ex}

For each $y_k\in H$, $k=1,2,\dots$, write
\[B(y_k):=\{\,x\in B^m(x_0,r_1):\, \forall  x_j\rightarrow x,\, x_j\in K_1,\, j=1,2,\dots,\,\Delta_{1*}u(x_j,y_k) {\textrm{ is convergent}}\,\}.\]
By (d),  $m_m(B^m(x_0,r_1)\setminus B(y_k))=0$ for all $k=1,2,\dots$. 
Write then 
\begin{equation*}B_1:=\bigcap_{k=1}^{+\infty}B(y_k)\setminus A.\end{equation*}  
Clearly $m_m(B^m(x_0,r_1)\setminus B_1)=0$. 
Similarly,  for each  $y\in B^n(y_0,r_2)$ write
\[B'(y):=\{\,x\in B^m(x_0,r_1):\, \forall  x_j\rightarrow x,\, x_j\in K_2,\, j=1,2,\dots,\,\exists  x_{j_l}, \, l=1,2,\dots ,
{\textrm{ s.t. }}\Delta_{1*}u(x_{j_l},y)\rightarrow \Delta_{1*}u(x,y)\,\}\]
and 
\begin{equation*}B_2:=\bigcap_{k=1}^{+\infty}B'(y_k),\,\, B=B_1\cap B_2,  {\textrm{ and }} G_2:=[(G\cap K_1\cap K_2)\setminus A]\cap (B_1\cap B_2).\end{equation*}  
By (e),  $m_m(B^m(x_0,r_1)\setminus B)=0$. One sees at once that $G_2\subset B$ is dense in  $B^m(x_0,r_1)$. 
To show the existence of the desired continuous extension, it is clearly sufficient to show that for each $(\tilde{x}_0,\tilde{y}_0)\in 
B\times B^n(y_0,r_2)$, the limit 
\begin{displaymath}\lim_{(x,y)\rightarrow (\tilde{x}_0,\tilde{y}_0),\, (x,y)\in G_2\times B^n(y_0,r_2)}\Delta _1u(x,y)\end{displaymath}
exists.  (This is of course standard, see e.g. [Di, (3.15.5), p.~54].) To see this, it is  sufficient to show that, 
for an arbitrary sequence $(x_j,y_j)\rightarrow (\tilde{x}_0,\tilde{y}_0)$, $(x_j,y_j)\in G_2\times B^n(y_0,r_2)$, $j=1,2, \dots$,  the limit
\begin{displaymath}\lim_{j\rightarrow +\infty }\Delta _1u(x_j,y_j)\end{displaymath}
exists. 
But this follows at once from the following facts:
\begin{enumerate}
\item[${\alpha )}$] The functions $\Delta _1u(x_j, \cdot )$, $j=1,2,\dots$, are nonnegative and harmonic in $B^n(y_0,r_2)$, by Step~3.
\item[${\beta )}$] For each $y_k\in H$, $k=1,2,\dots$,  the sequence $\Delta _1u(x_j,y_k)=\Delta _{1*}u(x_j,y_k)$ converges to 
$\Delta_{1*}u(\tilde{x}_0,y_k)=\Delta _1u(\tilde{x}_0,y_k)$ 
in $([0,+\infty ],q)$ as $j\rightarrow +\infty $. Moreover,  $\Delta _1u(\tilde{x}_0,y_k)\in {\mathbb{R}}$, (this follows from (d), (e) and from the fact 
that $x_j\in G_2$, $j=1,2,\dots$, and $\tilde{x}_0\in B$). 
\end{enumerate}
See  [Her71, Corollary~3, p.~6, and c), b), pp. 3, 2] (or [AG01, Lemma~1.5.6 and Lemma~1.5.10,  Theorem~1.5.8, pp.~16-17]),  and [V\" a71, Theorem~20.3,  p.~68]. 
By [AG01,  Theorem~1.5.8, p.~17]) one sees also the harmonicity of the functions
$\tilde{\Delta}_1u(x,\cdot ):B^n(y_0,r_2)\rightarrow {\mathbb{R}}, \,x\in B.$ 

\vspace{1ex}
\noindent{\textbf{Step~5.}} \emph{For each $x\in B^m(x_0,r_1)$ the  functions}
\begin{displaymath}B^n(y_0,r_2)\ni y\mapsto  \tilde{v}(x,y):=\int G_{B^m(x_0,r_1)}(x,z)\tilde{\Delta }_1u(z,y)dm_m(z)\in {\mathbb{R}}\end{displaymath}
\emph{and} 
\begin{displaymath}B^n(y_0,r_2)\ni y\mapsto  \tilde{h}(x,y):=u(x,y)+\tilde{v}(x,y)\in {\mathbb{R}}\end{displaymath}
\emph{are harmonic. Above and below $G_{B^m(x_0,r_1)}(x,z)$ is the Green function of the ball $B^m(x_0,r_1)$, with $x$ as a pole.}
 
\vspace{1ex}

Using  Fubini's Theorem one sees easily that for each $x\in B^m(x_0,r_1)$ the function 
 $\tilde{v}(x, \cdot )$ satisfies the mean value equality. To see that $\tilde{v}(x, \cdot )$ is harmonic, it is sufficient to show that 
$\tilde{v}(x, \cdot )\in {\mathcal{L}}^1_{\textrm{loc}}(B^n(y_0,r_2))$.   Using just   Fatou's Lemma, one sees that $\tilde{v}(x, \cdot )$ is 
lower semicontinuous, hence superharmonic.
Therefore either $\tilde{v}(x,\cdot )\equiv +\infty $ or else $\tilde{v}(x, \cdot )\in {\mathcal{L}}^1_{\textrm{loc}}(B^n(y_0,r_2))$. The following argument shows 
that the former alternative cannot 
occur. Indeed,   for each $x\in B$ and 
for each $y_k\in H$, $k=1,2,\dots$,
\[\tilde{\Delta}_{1}u(x,y_k)=\lim_{x'\rightarrow x,\,x'\in G_2}\Delta_{1}u(x',y_k)=\lim_{x'\rightarrow x,\,x'\in G_2}
\Delta_{1*}u(x',y_k)=\Delta_{1*}u(x,y_k)=\Delta_{1}u(x,y_k)\in {\mathbb{R}}.\]
 Hence $\tilde{v}(x,y_k)\in {\mathbb{R}}$  for each $x\in B^m(x_0,r_1)$ and $y_k\in H$, $k=1,2,\dots$.  Therefore, for each $x\in B^m(x_0,r_1)$ 
also the function $\tilde{h}(x,\cdot )=u(x,\cdot )+\tilde v(x,\cdot )$ is harmonic. 

\vspace{1ex}

\noindent{\textbf{Step~6.}} \emph{For each $y\in B^n(y_0,r_2)$ the  function}
\begin{displaymath}B^m(x_0,r_1)\ni x\mapsto  \tilde{h}(x,y):=u(x,y)+\tilde{v}(x,y)\in {\mathbb{R}}\end{displaymath}
\emph{is harmonic.}

\vspace{1ex}

With the aid of the version of Riesz's Decomposition Theorem, given in  [Ru50, 1.3, Theorem~II, p.~279, and p.~278, too] 
(see also [Sh56, Theorem~1, p.~499]), for  each $y\in B^n(y_0,r_2)$ one can write
\begin{displaymath}u(x,y)=h(x,y)-v(x,y)\end{displaymath}
where
\begin{displaymath}v(x,y):=\int G_{B^m(x_0,r_1)}(x,z)\Delta _1u(z,y)dm_m(z)\end{displaymath}
and $h(\cdot ,y)$ is the least harmonic majorant of $u(\cdot ,y)\mid B^m(x_0,r_1)$. Here  $v(\cdot ,y)$ is continuous and superharmonic in $B^m(x_0,r_1)$. 

As shown above in Step~5, $\Delta _1u(x,y_k)=\tilde{\Delta }_1u(x,y_k)$
for each  $x\in B$ and each  $y_k\in H$, $k=1,2,\dots$, and thus $v(\cdot ,y_k)=\tilde{v}(\cdot ,y_k)$ for each $y_k\in H$, $k=1,2,\dots $. 
Therefore  also $\tilde{h}(\cdot, y_k)=h(\cdot ,y_k)$ 
is harmonic for each $y_k\in H$, $k=1,2,\dots$.  To see that $\tilde{h}(\cdot ,y)$ is harmonic also for  $y\in B^n(y_0,r_2)\setminus H$, 
proceed in the following way.  
Let $B'(y)$ be the subset of $B^m(x_0,r_1)$, assumed in (e), see Step~4 above.
Write $G_3:=[(G\cap K_1\cap K_2)\setminus (A\cup A(y))]\cap (B\cap B'(y))$ 
and $B':=(B\cap B'(y))\setminus A(y)$. It is easy to see that  $G_3$ is dense in $B^m(x_0,r_1)$ and that $m_m(B^m(x_0,r_1)\setminus B')=0$. 
Take $x\in B'$ arbitrarily. Choose a sequence 
$x_j\rightarrow x$, $x_j\in G_3$, $j=1,2, \dots$.   By assumption~(e), there is a subsequence $x_{j_l}$ such that
\[\Delta_{1}u(x_{j_l},y)=\Delta_{1*}u(x_{j_l},y)\rightarrow \Delta _{1*}u(x,y)=\Delta_{1}u(x,y)\]
 as $l\rightarrow +\infty $. On the other hand, since $x_{j_l}\in G_3$, $l=1,2,\dots $,
\[\Delta_{1}u(x_{j_l},y)=\tilde{\Delta}_{1}u(x_{j_l},y)\rightarrow \tilde{\Delta}_{1}u(x,y),\] we 
have $\Delta _{1}u(x,y)=\tilde{\Delta}_{1}u(x,y)$. Since  $m_m(B^m(x_0,r_1)\setminus B')=0$, we have
$v(x,y)=\tilde{v}(x,y)$, and thus  $\tilde {h}(x,y)=h(x,y)$, for all $x\in B^m(x_0,r_1)$ and all
$y\in B^n(y_0,r_2)$. Therefore $\tilde{h}(\cdot ,y)$ is harmonic for all $y\in B^n(y_0,r_2)$. 

\vspace{1ex}

\noindent{\textbf{Step~7.}} \emph{The use of the results of  Lelong and of Avanissian.}

\vspace{1ex}

By Steps~5 and 6 we know that $h(\cdot ,\cdot )=\tilde{h}(\cdot ,\cdot )$ is separately harmonic in $B^m(x_0,r_1)\times B^n(y_0,r_2)$.
By Lelong's result [Le61, Théorème~11, p.~554] or [Av67, Théorème~1, pp.~4-5] (for a more general result, 
see [Si69, Theorem~7.1, p.~166, and Corollary, p.~145], see also [Im90, Theorem, p.~9]), $\tilde{h}(\cdot ,\cdot )$ is harmonic and thus 
locally bounded above in $B^m(x_0,r_1)\times B^n(y_0,r_2)$. Therefore also $u(\cdot ,\cdot )$ is locally bounded above in $B^m(x_0,r_1)\times B^n(y_0,r_2)$.
But then it follows from Avanissian's result [Av61, Théorème~9, p. 140] (or [Ar66, Theorem~1, p.~622] or [Ri89,Theorem~1, p. 69]) that $u(\cdot ,\cdot )$ is subharmonic on $B^m(x_0,r_1)\times B^n(y_0,r_2)$.
\hfill \qed
 
\vspace{2ex}
 
\noindent{\textbf{5.3.}} Another variant of the above result is the following, where the assumption (e) is replaced with a certain ``continuity''
 condition of $v(\cdot ,\cdot )$ in the second variable. Our result is an improved version of the result [Ri07$_2$, Theorem~2, pp.~442-443]: 
Again the conditions (d) and (e) are essentially milder and consequently the proof contains new ingredients.

\vspace{1ex}

\noindent{\textbf{Theorem~5.}}  \emph{ Let $\Omega $ be a domain in ${\mathbb{R}}^{m+n}$, \mbox{$m,n\geq 2$.} 
Let $u:\, \Omega \rightarrow 
{\mathbb{R}}$ be such that for each $(x',y')\in \Omega $ there is  $(x_0,y_0)\in \Omega $ and $r_1>0$, $r_2>0$  such that 
$(x',y')\in B^m(x_0,r_1)\times B^n(y_0,r_2)\subset \overline{B^m(x_0,r_1)\times B^n(y_0,r_2)}\subset \Omega $ and  such that the following 
conditions are satisfied:} 
 \begin{itemize}
\item[(a)] \emph{For each $y\in \overline{B^n(y_0,r_2)}$ the function} 
\[\overline{B^m(x_0,r_1)}\ni x\mapsto u(x,y)\in {\mathbb{R}}\]
\emph{is continuous, and subharmonic in $B^m(x_0,r_1)$.}
 \item[(b)] \emph{For each $x\in \overline{B^m(x_0,r_1)}$ the function} 
\[\overline{B^n(y_0,r_2)}\ni y\mapsto u(x,y)\in {\mathbb{R}}\]
\emph{is continuous, and harmonic in $B^n(y_0,r_2)$.}
\item[(c)] \emph{For each $y\in B^n(y_0,r_2)$ one has $\Delta _{1*}u(x,y)<+\infty $ for every $x\in B^m(x_0,r_1)$, possibly with the exception of 
a polar set in $B^m(x_0,r_1)$.}
\item[(d)] \emph{There is a set $H\subset B^n(y_0,r_2)$, dense in $B^n(y_0,r_2)$, and an open  set $K_1\subset B^m(x_0,r_1)$, dense in $B^m(x_0,r_1)$,
such that} 

\vspace{1ex}

\begin{itemize}
\item[(d1)]  \emph{for each $y\in H$, for almost every $x\in B^m(x_0,r_1)$ and for any sequence $x_j\rightarrow x$, $x_j\in K_1$, $j=1,2,\dots$, one has 
\[\Delta_{1*}u(x_j,y)\rightarrow \Delta_{1*}u(x,y)\]
in $([0,+\infty ],q)$ as $j\rightarrow +\infty$,}
\item[(d2)] \emph{for some $y\in H$ one has}
\[\sup_{x\in K_1}\Delta_{1*}u(x,y)<+\infty .\]
\end{itemize}

\vspace{1ex}

\item[(e)] \emph{There is a set $K_2\subset B^m(x_0,r_1)$, dense in $B^m(x_0,r_1)$, such that for each $x\in K_2$ and 
for each $y\in B^n(y_0,r_2)$ there is a sequence $y_k\rightarrow y$, $y_k\in H$, $k=1,2,\dots$,  such that  
\[\int G_{B^m(x_0,r_1)}(x,z)\, \Delta _1u(z,y_{k})\, dm_m(z)\rightarrow \int G_{B^m(x_0,r_1)}(x,z)\, \Delta _1u(z,y)\, dm_m(z)\]
as $k\rightarrow +\infty $.}
\end{itemize}
\emph{Then} $u$ 
\emph{is subharmonic in $\Omega $}.

\vspace{1ex}
\noindent{\emph{Proof}}. The proof differs from the proof of Theorem~4 above only in Steps~4 and 6, in the following way. 

\vspace{1ex}

\noindent{\textbf{Step~4.}} \emph{The function  $\Delta_{1}u(\cdot ,\cdot )\mid G'_2\times B^n(y_0,r_2)$ 
has a continuous and bounded extension $\tilde{\Delta}_{1}u(\cdot ,\cdot ): B_1\times B^n(y_0,r_2)\rightarrow {\mathbb{R}}$ 
(for the definition of $B_1$, see above the proof of Theorem~4,  for the definition of $G'_2$ see below).  Moreover, the functions 
$\tilde{\Delta}_1u(x,\cdot )$, $x\in B_1$, are nonnegative and harmonic  in $B^n(y_0,r_2)$, and  one may suppose that they are uniformly bounded  
in $B^n(y_0,r_2)$ by a fixed constant $M$.} 
 
\vspace{1ex}

Write $G'_2:=G\cap K_1\cap B_1$. Since $m_m(B^m(x_0,r_1)\setminus B_1)=0$, one sees at once that $G'_2$ is dense in  $B^m(x_0,r_1)$. 
To show that $\Delta_1u(\cdot ,\cdot )\mid G'_2\times B^n(y_0,r_2)$ has a continuous extension
$\tilde{\Delta}_1u(\cdot ,\cdot ): \, B_1\times B^n(y_0,r_2)\rightarrow {\mathbb{R}}$, just proceed as above in the proof of Theorem~4, Step~4.  

Also the harmonicity of the functions
$\tilde{\Delta}_1u(x,\cdot ), \,x\in B_1$,
follows as above. In addition, because of  assumption (d2),  one sees that now the family
\begin{displaymath}\tilde{\Delta}_1u(x,\cdot )=\Delta_1u(x,\cdot ), \,x\in G'_2,\end{displaymath}
and thus also the family
\begin{displaymath}\tilde{\Delta}_1u(x,\cdot ), \,x\in B_1,\end{displaymath}
of nonnegative  harmonic functions in $B^n(y_0,r_2)$, is uniformly locally bounded in $B^n(y_0,r_2)$. As a matter of fact (see the beginning of the proof of Theorem~4), 
we may suppose 
that these families are  even uniformly  bounded  in $B^n(y_0,r_2)$, by a fixed constant $M$, say. Therefore also  $\tilde{\Delta}_{1}u(\cdot , \cdot )$ is 
bounded in $B_1\times B^n(y_0,r_2)$ by this same $M$. 

\vspace{1ex}

\noindent{\textbf{Step~6.}} \emph{For each $y\in B^n(y_0,r_2)$ the  function}
\begin{displaymath}B^m(x_0,r_1)\ni x\mapsto  \tilde{h}(x,y):=u(x,y)+\tilde{v}(x,y)\in {\mathbb{R}}\end{displaymath}
\emph{is harmonic.}

\vspace{1ex}

With the aid of the cited version of Riesz's Decomposition Theorem, one can again,  for  each $y\in B^n(y_0,r_2)$, write
\begin{displaymath}u(x,y)=h(x,y)-v(x,y)\end{displaymath}
where $v(\cdot ,y)$ is continuous and superharmonic in $B^m(x_0,r_1)$ 
and $h(\cdot ,y)$ is the least harmonic majorant of $u(\cdot ,y)\mid B^m(x_0,r_1)$. 
As in  the proof of Theorem~4, Step~6, one sees that for all $y_k\in H$, $k=1,2,\dots$, the function $\tilde{h}(\cdot ,y_k)=h(\cdot ,y_k)$ is harmonic in $B^m(x_0,r_1)$.   
To see that $\tilde{h}(\cdot ,y)$ is harmonic also for 
$y\in B^n(y_0,r_2)\setminus H$, just use 
the present assumption (e),  in the following way.  Choose $y\in B^n(y_0,r_2)$ and $x\in K_2$ arbitrarily. Take then  a sequence 
$y'_k\rightarrow y$, $y'_k\in H$, $k=1,2,\dots$.
Then the functions $\tilde{h}(\cdot ,y'_k)=h(\cdot ,y'_k)$, $k=1,2,\dots$,
are harmonic in $B^m(x_0,r_1)$.  By the assumption (e), $v(x,y'_k)\rightarrow v(x,y)$ as $k\rightarrow +\infty $. 
On the other hand, since $\tilde{v}(x,\cdot )$ is harmonic, $\tilde{v}(x,y'_k)\rightarrow \tilde{v}(x,y)$ as $k\rightarrow +\infty $. Since $v(x,y'_k)
=\tilde{v}(x,y'_k)$, $k=1,2,\dots$, we see that $v(x,y)=\tilde{v}(x,y)$ for each $x\in K_2$. Using then the facts that $v(\cdot ,y)$ is continuous, that, by          
[Hel69, Theorem~6.22, p.~119], say, also $\tilde{v}(\cdot ,y)$ is continuous (recall that $\tilde{\Delta}_1u(\cdot ,\cdot )$ is continuous and bounded in 
$B_1\times B^n(y_0,r_2)$, and that $m_m(B^m(x_0,r_1)\setminus B_1)=0$, see Step~4 above), and that $K_2$ is dense in $B^m(x_0,r_2)$, 
we see that $v(\cdot ,y)=\tilde{v}(\cdot ,y)$. Therefore  for all $x\in B^m(x_0,r_1)$,
\begin{displaymath}\tilde{h}(x,y)=u(x,y)+\tilde v(x,y)=u(x,y)+v(x,y)=h(x,y),\end{displaymath} 
and the harmonicity of $\tilde{h}(\cdot ,y)$ in $B^m(x_0,r_1)$ follows. 

Now we know that $\tilde{h}(\cdot ,\cdot )$ is separately harmonic in $B^m(x_0,r_1)\times B^n(y_0,r_2)$. The rest of the proof goes then as in Step~7 of 
the proof of Theorem~4 above.\hfill \qed

\vspace{2ex}

\noindent{\textbf{5.4.}} The assumptions of Theorems~4 and 5 above, especially the (e)-assumptions,  are undoubtedly still somewhat technical. 
However, replacing  Ko\l odziej's and Thorbi\"ornson's ${\mathcal{C}}^2$ assumption of the functions $u(\cdot ,y)$ by the (spherical) continuity 
requirement of the 
generalized Laplacians $\Delta _1u(\cdot ,y)$, we  obtain the following  concise corollaries to Theorem~4. Observe that generalized Laplacians of subharmonic 
functions may indeed assume also the value $+\infty$, see  Example~2 in {\textbf{5.1}} above.

\vspace{1ex}

\noindent{\textbf{Corollary~1.}} \emph{ Let $\Omega $ be a domain in ${\mathbb{R}}^{m+n}$, \mbox{$m,n\geq 2$.} 
Let $u:\, \Omega \rightarrow 
{\mathbb{R}}$ be such that} 
 \begin{itemize}
\item[(a)] \emph{for each $y\in {\mathbb{R}}^n$ the function} 
\[\Omega (y)\ni x\mapsto u(x,y)\in {\mathbb{R}}\]
\emph{is continuous and subharmonic,}
 \item[(b)] \emph{for each $x\in {\mathbb{R}}^m$ the function} 
\[\Omega (x)\ni y\mapsto u(x,y)\in {\mathbb{R}}\]
\emph{is harmonic,}
\item[(c)] \emph{for each $y\in {\mathbb{R}}^n$ the function} 
\[\Omega (y)\ni x\mapsto \Delta _1u(x,y)\in [0,+\infty ]\]
\emph{is defined,  continuous (with respect to the spherical metric), and finite for all $x$, except at most of a polar set $E(y)$ in $\Omega (y)$.} 
\end{itemize}
\emph{Then} $u$ 
\emph{is subharmonic in $\Omega $}.

\vspace{1ex}

\noindent{\textbf{Corollary~2.}} ([Ri07$_2$, Corollary, p.~444])  \emph{ Let $\Omega $ be a domain in ${\mathbb{R}}^{m+n}$, \mbox{$m,n\geq 2$.} 
Let $u:\, \Omega \rightarrow 
{\mathbb{R}}$ be such that} 
 \begin{itemize}
\item[(a)] \emph{for each $y\in {\mathbb{R}}^n$ the function} 
\[\Omega (y)\ni x\mapsto u(x,y)\in {\mathbb{R}}\]
\emph{is continuous and subharmonic,}
 \item[(b)] \emph{for each $x\in {\mathbb{R}}^m$ the function} 
\[\Omega (x)\ni y\mapsto u(x,y)\in {\mathbb{R}}\]
\emph{is harmonic,}
\item[(c)] \emph{for each  $y\in {\mathbb{R}}^n$ the function} 
\[\Omega (y)\ni x\mapsto \Delta _1u(x,y)\in {\mathbb{R}}\]
\emph{is defined  and continuous.}
\end{itemize}
\emph{Then} $u$ 
\emph{is subharmonic in $\Omega $}.

\vspace{3ex}

\flushleft{\noindent\textbf{References}}

\vspace{1ex}
\begin{flushleft}
\begin{enumerate}
\item[{[ABR01]}] Axel, S., Bourdon, P., Ramey, W.  ``Harmonic Function Theory'',  
Springer-Verlag, New-York, 2001 (Second Edition).
\item[{[AG93]}] Armitage, D.H., Gardiner, S.J.
 ``Conditions for separately subharmonic functions to be subharmonic'',  
Potential Anal.,
{\textbf{2}} (1993), 255--261.
\item[{[AG01]}] Armitage, D.H., Gardiner, S.J.
 ``Classical Potential Theory'',  
Springer-Verlag, London, 2001.
\item[{[Ar66]}] Arsove, M.G. ``On subharmonicity of doubly subharmonic functions'',  
Proc. Amer. Math. Soc.,
{\textbf{17}} (1966), 622--626.
\item[{[Av61]}] Avanissian, V.
 ``Fonctions plurisousharmoniques et fonctions doublement sousharmoniques'',  
Ann. Sci. École  Norm. Sup., {\textbf{78}} (1961), 101--161.
\item[{[Av67]}] Avanissian, V.
 ``Sur l'harmonicité des fonctions séparément harmoniques'', in: 
Séminaire de Probabilités (Univ. Strasbourg, Février 1967), {\textbf{1}} (1966/1967), pp.~101--161, Springer-Verlag, Berlin, 1967.
\item[{[Br38]}] Brelot, ~M. 
 ``Sur le potentiel et les suites de fonctions sousharmoniques'',  
C.R. Acad. Sci., {\textbf{207}} (1938), 836--838.
\item[{[Br69]}] Brelot, M.
 ``Éléments de la Théorie Classique du Potentiel'',  
Centre de Documentation Universitaire, Paris, 1969 (Third Edition).
\item[{[CZ61]}] Calderon, A.P., Zygmund, A. ``Local properties of solutions of elliptic partial differential equations'', 
Studia Math., {\textbf{20}} (1961), 171--225.
\item[{[CS93]}] Cegrell, U.,  Sadullaev, A. ``Separately subharmonic functions'', Uzbek. Math. J., {\textbf{1}} (1993), 78--83.
\item[{[Di60]}] Dieudonné, J.
 ``Foundations of Modern Analysis'',  
Academic Press, New York, 1960.
\item[{[Do57]}] Domar, Y. ``On the existence of a largest subharmonic minorant of a given function'', 
Arkiv f\"or matematik, {\textbf{3}}, nr. 39 (1957), 429--440.
\item[{[Do88]}] Domar, Y. ``Uniform boundedness in families related to subharmonic functions'', 
J. London Math. Soc. (2), {\textbf{38}} (1988), 485--491.
\item[{[Hel69]}] Helms, L.L. ``Introduction to Potential Theory'', 
Wiley-Interscience, New York, 1969.
\item[{[Her71]}] Hervé, M. ``Analytic and Plurisubharmonic Functions in Finite and Infinite Dimensional Spaces'', 
Lecture Notes in Mathematics 198, Springer-Verlag, Berlin, 1971.  
\item[{[Im90]}] Imomkulov, S.A. ``Separately subharmonic functions'' (in Russian),  
Dokl. USSR, {\textbf{2}} (1990), 8--10.
\item[{[KT96]}] Ko\l odziej, S., Thorbi$\ddot {\textrm{o}}$rnson, J. ``Separately harmonic and subharmonic functions'',  
Potential Anal.,
{\textbf{5}} (1996), 463--466.
\item[{[Le45]}] Lelong, P. ``Les fonctions plurisousharmoniques'',  
Ann. Sci. École Norm. Sup.,
{\textbf{62}} (1945), 301--338.
\item[{[Le61]}] Lelong, P. ``Fonctions plurisousharmoniques et fonctions analytiques de variables réelles'',  
Ann. Inst. Fourier, Grenoble, 
{\textbf{11}} (1961), 515--562.
\item[{[Le69]}] Lelong, P. ``Plurisubharmonic Functions and Positive Differential Forms'',  
Gordon and Breach, London, 1969.
\item[{[LL01]}] Lieb, E.H., Loss, M. ``Analysis'',
 Graduate Studies in Mathematics, 14, American Mathematical Society, Providence, Rhode Island, 2001.
\item[{[Ma95]}] Mattila,~P. ``Geometry of Sets and Measures in Euclidean Spaces'',
Cambridge studies in advanced mathematics, 44, Cambridge University Press, Cambridge,  1995.
\item[{[Mi96]}] Mizuta,~Y. ``Potential Theory in Euclidean Spaces'',
 Gaguto International Series, Mathematical Sciences and Applications, 6, Gakk$\bar{{\textrm{o}}}$tosho Co., Tokyo, 1996.
\item[{[Pa94]}] Pavlovi\'c, M. ``On subharmonic behavior and oscillation of functions on balls in ${\mathbb{R}}^n$'', 
Publ. Inst. Math. (Beograd), {\textbf{55 (69)}} (1994), 18--22.
\item[{[PR08]}] Pavlovi\'c, M., Riihentaus, J. ``Classes of quasi-nearly subharmonic functions'', Potential Anal., {\textbf{29}} (2008), 89--104.
\item[{[Pl70]}] du Plessis, N. ``An Introduction to Potential Theory'', Oliver \& Boyd, Edinburgh, 1970.
\item[{[Ra37]}] Radó, T. ``Subharmonic Functions'',
 Springer-Verlag, Berlin, 1937.
\item[{[Ri84]}] Riihentaus, J. ``On the extension of separately hyperharmonic functions and H$^p$-functions'',
Michigan Math. J., {\textbf{31}} (1984),  \mbox{99--112.}
\item[{[Ri89]}] Riihentaus, J. ``On a theorem of Avanissian--Arsove'',
Expo.  Math., {\textbf{7}} (1989),  \mbox{69--72.}
\item[{[Ri00]}] Riihentaus, J.  ``Subharmonic functions: non-tangential and 
tangential boundary
behavior'' in:   Function Spaces, Differential Operators and Nonlinear Analysis (FSDONA'99), Proceedings of the Sy\"ote Conference 1999, 
  Mustonen, V.,  R\'akosnik, J. (eds.), 
 Math. Inst., Czech Acad. Science,  Praha, 2000, \mbox{pp. 229--238.} \mbox{(ISBN 80-85823-42-X)}
\item[{[Ri01]}] Riihentaus, J.  ``A generalized mean value inequality for subharmonic functions'', Expo. Math., {\textbf{19}} (2001),  \mbox{187-190}.
\item[{[Ri05]}] Riihentaus, J. ``An integrability condition and weighted boundary behavior of subharmonic and ${\mathcal{M}}$-subharmonic functions:
a survey'', Int. J. Diff. Eq. Appl., {\textbf{10}} (2005), 1--14.
\item[{[Ri06$_1$]}] Riihentaus, J.  ``A weighted boundary limit result for subharmonic functions'', Adv. Algebra and Analysis, {\textbf{1}} (2006), 
\mbox{27--38.}
\item[{[Ri06$_2$]}] Riihentaus, J. ``Separately subharmonic functions'', arXiv:math.AP/0610259 v3 13 Nov 2006. 
\item[{[Ri06$_3$]}] Riihentaus, J. ``Separately quasi-nearly subharmonic functions'', in: Complex Analysis and Potential Theory, Proceedings of the 
Conference Satellite to ICM~2006, Tahir Aliyev Azero$\breve{\textrm{g}}$lu, Promarz M. Tamrazov (eds.), Gebze Institute of Technology, Gebze, Turkey, 
September 8-14,
 2006, World Scientific, Singapore, 2007, pp.~156-165.
\item[{[Ri07$_1$]}] Riihentaus, J. ``On the subharmonicity of separately  subharmonic functions'', in: Proceedings of the 11th WSEAS International 
Conference on Applied Mathematics (MATH'07),  Dallas, Texas, USA, March 22-24, 2007, Kleanthis Psarris, Andrew D.~Jones (eds.), WSEAS, 2007,  
pp.~230-236. 
(IBSN 978-960-8457-60-7)
\item[{[Ri07$_2$]}] Riihentaus, J. ``On separately harmonic and subharmonic functions'', Int. J. Pure Appl. Math., {\textbf{35}}, no. {\textbf{4}} (2007), 
\mbox{435-446}.
\item[{[Ru50]}] Rudin, W. ``Integral representation of continuous functions'', Trans. Amer. Math. Soc., {\textbf{68}} (1950), 
278--286.
\item[{[Ru79]}] Rudin, W. ``Real and Complex Analysis'', 
Tata McGraw-Hill, New Delhi, 1979.
\item[{[Sa41]}] Saks, S. ``On the operators of Blaschke and Privaloff for subharmonic functions'', Rec. Math. (Mat. Sbornik), {\textbf{9 (51)}} (1941), 
451--456.
\item[{[Sh56]}] Shapiro, V.L.  ``Generalized laplacians'',  Amer. J. Math.,
{\textbf{78}} (1956), \mbox{497--508.}
\item[{[Sh71]}] Shapiro, V.L.  ``Removable sets for pointwise subharmonic functions'', Trans. Amer. Math. Soc.,
{\textbf{159}} (1971), \mbox{369--380.}
\item[{[Sh78]}] Shapiro, V.L.  ``Subharmonic functions and Hausdorff measure'', J. Diff. Eq., {\textbf{27}} (1978), \mbox{28--45.}
\item[{[Si69]}] Siciak, J.  ``Separately analytic functions and envelopes of holomorphy of some lower dimensional subsets of ${\mathbb{C}}^n$'', 
Ann. Polon. Math.,
{\textbf{22}} (1969), \mbox{145--171.}
\item[{[Sz33]}] Szpilrajn, E.  ``Remarques sur les fonctions sousharmoniques'',  Ann. Math., {\textbf{34}} (1933), \mbox{588--594.}
\item[{[V\" a71]}] V\"ais\"al\"a, J. ``Lectures on n-Dimensional Quasiconformal Mappings'', 
Lecture Notes in Mathematics 229, Springer-Verlag, Berlin, 1971.  
\item[{[Vu82]}] Vuorinen, M.  ``On the Harnack constant and the boundary behavior of Harnack functions'',  Ann. Acad. Sci. Fenn., Ser. A I, 
Math., {\textbf{7}} (1982), 259--277.
\item[{[Wi88]}] Wiegerinck, J.  ``Separately  subharmonic functions need not be subharmonic'',  
Proc.  Amer. Math.  Soc., {\textbf{104}} (1988), 770-771.
\item[{[WZ91]}] Wiegerinck, J., Zeinstra, R.  ``Separately subharmonic functions: when are they subharmonic'' in: Proceedings of Symposia in Pure 
Mathematics, vol. {\textbf{52}}, part {\textbf{1}},  Eric Bedford, John P. D'Angelo, Robert E.~Greene, Steven G.~Krantz (eds.), 
 Amer. Math. Soc., Providence, Rhode Island, 1991,  \mbox{pp. 245--249}.

\end{enumerate}
\end{flushleft}
\end{document}